\documentclass{amsart}

\usepackage{amsmath,amsthm,amsfonts,amssymb,bm,graphicx, mathrsfs}

\usepackage
{hyperref}
\hypersetup{colorlinks=true,citecolor=blue,linkcolor=blue,urlcolor=blue,
pdfstartview=FitH }

\theoremstyle{plain}
\newtheorem{theorem}{Theorem}

\newtheorem{lemma}{Lemma}

\numberwithin{equation}{section}

\theoremstyle{definition}

\renewcommand{\geq}{\geqslant}
\renewcommand{\leq}{\leqslant}

\newcommand{\changed}[1]{{\color{black} #1}}

\usepackage{calc}
\newsavebox\CBox
\newcommand\hcancel[2][0.5pt]{%
  \changed{\ifmmode\sbox\CBox{$#2$}\else\sbox\CBox{#2}\fi%
  \makebox[0pt][l]{\usebox\CBox}%
  \rule[0.5\ht\CBox-#1/2]{\wd\CBox}{#1}}}

\makeatletter
\DeclareRobustCommand\widecheck[1]{{\mathpalette\@widecheck{#1}}}
\def\@widecheck#1#2{%
    \setbox\z@\hbox{\m@th$#1#2$}%
    \setbox\tw@\hbox{\m@th$#1%
       \widehat{%
          \vrule\@width\z@\@height\ht\z@
          \vrule\@height\z@\@width\wd\z@}$}%
    \dp\tw@-\ht\z@
    \@tempdima\ht\z@ \advance\@tempdima2\ht\tw@ \divide\@tempdima\thr@@
    \setbox\tw@\hbox{%
       \raise\@tempdima\hbox{\scalebox{1}[-1]{\lower\@tempdima\box
\tw@}}}%
    {\ooalign{\box\tw@ \cr \box\z@}}}
\makeatother

\begin{document}

\author{Valentin Blomer}

\address{Mathematisches Institut, Bunsenstr. 3-5, 37073 G\"ottingen, Germany} \email{vblomer@math.uni-goettingen.de}

 \author{Rizwanur Khan}
 \address{Science Program, Texas A\&M University at Qatar, PO Box 23874, Doha, Qatar}
\email{rizwanur.khan@qatar.tamu.edu}

 \title{Uniform subconvexity and symmetry breaking reciprocity }

\thanks{First author supported in part by DFG grant BL 915/2-2.}
 

\begin{abstract}  A non-symmetric reciprocity formula is established that expresses  the fourth moment of automorphic  $L$-functions of level $q$ and primitive central character  twisted by the $\ell$-th Hecke eigenvalue as a twisted mixed   moment of automorphic  $L$-functions of level $\ell$ and trivial central character.  As an application, uniform  subconvexity bounds for $L$-functions in the level and the eigenvalue aspect are derived. 
\end{abstract}


\setcounter{tocdepth}{2}  \maketitle 

\maketitle

 \section{Introduction} 
 
 \subsection{Subconvexity} Although the subconvexity problem for automorphic $L$-functions on ${\rm GL}(2)$ is by and large well-understood and -- without an explicit exponent --  solved in all aspects and over all  number fields simultaneously \cite{MV}, some aspects are harder than others. For a  Maa{\ss} form $f$ of level $q$ and (necessarily even) primitive central character $\chi$,  subconvexity in the $q$-aspect was only solved in 2002  in a 90-page paper \cite{DFI} facing enormous technical and structural difficulties and saving only a power $q^{1/23041}$ from the trivial bound $q^{1/4}$. This is rather surprising, as much simpler and numerically much stronger results were already available in the case of trivial central character (as well as for holomorphic forms and  subconvexity in other aspects). The proof in \cite{DFI} was simplified, generalized and improved in \cite{BHM} and became an important ingredient in a cubic analogue of Duke's equidistribution result \cite{ELMV}. Yet the proof remained complicated and the saving $q^{1/1889}$  relatively microscopic, supporting the empirical fact that leaving the realm of trivial central character may offer all sorts of difficult phenomena (see \cite{IL} for a particularly striking result). In this paper we   take a fresh look at the hardest of the ${\rm GL}(2)$ subconvexity problems and prove the following \emph{uniform} subconvexity result. Here and henceforth, $\vartheta_0 \leq 7/64$ denotes an admissible exponent for the Selberg eigenvalue conjecture (i.e.\ the Ramanujan conjecture only at the place $\infty$), and $\vartheta \leq 7/64$ denotes an admissible exponent for the ${\rm GL}(2)$ Ramanujan eigenvalue conjecture for all places. 
 
 \begin{theorem}\label{thm1} Let $q$ be a  prime, $\tau \in \Bbb{R}$, $t\in \Bbb{R} \cup [-i\vartheta_0, i\vartheta_0]$, and write $T = 1+|t|$.  Let $f$ be a Hecke-Maa{\ss} cusp form of eigenvalue $\lambda = 1/4 + t^2$, level $q$ and primitive central character. Then
 $$L(1/2+i\tau, f) \ll_{\tau, \varepsilon} (qT)^{ \varepsilon}  \left(q^{\frac{1}{4} - \frac{1}{128}} T^{\frac{1}{2} - \frac{1 - 2\vartheta_0}{20}} + q^{\frac{1}{8}} T^{\frac{1}{2}} \right)$$
 for any $\varepsilon > 0$. 
 \end{theorem}

The same result with essentially the same proof holds for holomorphic forms (where the weight $k$ plays the role of $T$), but as this case is known to be easier \cite{DFI3}, we focus here on the Maa{\ss} case. 

We regard $q$ as the main parameter, and our result showcases a much superior 
   saving in the $q$-aspect compared to both \cite{DFI} and \cite{BHM}. We did not make an effort to optimize the $T$-aspect, nevertheless  our result is also quite strong in the $T$-aspect and hence  simultaneously in    \emph{all} defining parameters of the automorphic form $f$. In fact, we obtain uniform subconvexity in the analytic conductor $\mathcal{C}(f) = qT^2$  of $f$ unless $q$ is (essentially) fixed and only the archimedean parameter $T$ grows. This case, however, has been well-understood for a long time; for level 1 there are Weyl-type savings in $T$ in the literature  \cite{Iv} (see also \cite{Iw1} for different method with a slightly weaker exponent) that could be generalized straightforwardly to level $q$ with polynomial dependence in $q$. 
       
 
More important than the exact shape of the bound is the new set of methods that are substantially different from all previous approaches to the subconvexity problem with nebentypus. The proof of Theorem \ref{thm1} is based on a new spectral reciprocity formula that we will describe in the next subsection. At this point we remark that we start rather classically with an amplified fourth moment, but we emphasize two important differences compared to previous approaches. On the technical side, we never use approximate functional equations, but work with complete Dirichlet series until the very end and let instead analytic continuation do the work for us. This approach is much more flexible and one of the reasons why a uniform saving comes almost for free\footnote{Cf.\   \cite[p.\ 702]{BHM}: ``\emph{We believe that a line of attack that totally dispenses with approximate functional equations (as in \cite{By}) would give a cleaner and simpler proof, but we have not yet succeeded in completing this project.}''}. We do make use, however, of the functional equation right from the start by considering a fourth moment of terms
$$L(s , f)^2 L(1-s , \bar{f}) L(w, \bar{f})$$
for $\Re s,   \Re w$ sufficiently large. An application of the functional equation brings us in the region of absolute convergence (but more importantly introduces the root number) and at the end we continue analytically to $s= \frac{1}{2} +  i\tau$, $w = \frac{1}{2} - i\tau$. A more detailed description of this device will be given later. 
On the structural side, we use higher rank methods, and in particular a version of the ${\rm GL}(3)$ Voronoi summation formula, even though Theorem \ref{thm1} can be stated purely in terms of ${\rm GL}(2)$. 

\subsection{Reciprocity} Reciprocity formulae for $L$-functions express moments of $L$-functions in terms of different moments of $L$-functions. Maybe the first such identity goes back to Kuznetsov and Motohashi  \cite{Mo} featuring a fourth moment of level 1 automorphic forms. Motohashi's famous formula relating the fourth moment of the Riemann zeta function \cite{Mo1} can also be seen as a reciprocity formula of families of $L$-functions. An example for Dirichlet $L$-functions was found by Conrey and refined by Young \cite{Yo} and Bettin \cite{Be}. Various automorphic reciprocity formulae were recently established in \cite{BLM, BK, AK}; all of the latter consider products of $L$-functions of total degree 8 and are symmetric in the sense that the families of $L$-functions on both sides are of the same type. Here we develop a \emph{symmetry breaking} reciprocity formula roughly of the following type
\begin{equation}\label{proto}
\sum_{f \in \mathcal{B}(q, \chi)} |L(1/2, f )|^4   \lambda_f(\ell) \rightsquigarrow \sum_{f \in \mathcal{B}(\ell, \text{triv})} L(1/2, f)^3 L(1/2, f \times \bar{\chi}) 
\end{equation}
where the left hand side runs through Maa{\ss} forms of level $q$ and primitive central character $\chi$, while the right hand runs through Maa{\ss} forms of level $\ell$ and trivial central character. In contrast to the formula in \cite{BK}, the right hand side is not twisted by $\lambda_f(q)$, but the twist is rather ``inside the $L$-function'' by the central character. (As an aside,   if $\chi$ is quadratic, the factor $L(1/2, f \times \bar{\chi})$ can be seen as the square of a  metaplectic Fourier coefficient of index $q$.) Our formula also makes precise the archimedean weight functions that are suppressed in the above notation. Careful estimations along with the best (hybrid) subconvexity bounds for $L(1/2, f \times \bar{\chi})$ produce the bounds in Theorem \ref{thm1}. The connection is roughly as follows: if 
\begin{equation}\label{twist}
L(1/2, f \times \bar{\chi}) \ll q^{\alpha} N^{\beta} (1 + |t_f|)^{O(1)}
\end{equation}
for $f \in \mathcal{B}(N, \text{triv})$ 
with $0 < \alpha < 1/2$ and $\beta > 0$, then
$$L(1/2, f) \ll q^{\frac{1}{4} - \frac{1 - 2\alpha}{8(3 + 4\beta)}+\varepsilon}  (1 + |t_f|)^{\frac{1}{2} - \frac{1-2\vartheta_0}{4(3+4\beta)}+\varepsilon}$$
for $f \in \mathcal{B}(q, \chi)$; at least this is the limit of the method. It is another instance where a subconvexity bound for an ``easier'' family of $L$-function becomes instrumental for a subconvexity bound of a ``more complicated'' family, cf.\ \cite{MV, Ne} as well as the ancestor of this paper, \cite{BHM}, where also a bound of type \eqref{twist} was used.

That a formula of type \eqref{proto} should exist, can be seen by a heuristic back-of-an-envelope computation, cf.\ Section \ref{heu}. 
This computation assumes that ``everything is coprime to everything'', which should morally be not too far from reality, but making this precise seems to be incredibly difficult, because each step introduces a whole alphabet of auxiliary variables. Whoever has tried to find an elementary  proof of Selberg's formula 
$$S(a, b, c) = \sum_{d \mid (a, b, c)} d S(1, ab/d^{2}, c/d)$$ 
for Kloosterman sums (i.e.\ without interpreting it via spectral theory as a statement of the ${\rm GL}(2)$ Hecke algebra)
 has seen a toy case of a situation where a clean and seemingly elementary identity of exponential sums appears to resist straightforward attempts\footnote{See e.g.\ \cite{An, HK} an the references therein for short elementary proofs.}.  
 As in \cite{BK} the key to a clean approach is to use higher rank tools and to replace three Poisson summations by the ${\rm GL}(3)$ Voronoi summation. This packages the Hecke combinatorics nicely and cleanly, and it has the additional advantage that it works also in the cuspidal case, and not only in the Eisenstein case. The new problem, however, is that we would have to apply a Voronoi summation formula to an Eisenstein series of level $q$, whose Hecke eigenvalues at a prime $p \nmid q $ are given by $(\textbf{1} \ast \textbf{1} \ast \bar{\chi})(p),$ 
and again the combinatorial issues become extremely convoluted. 

 Therefore we apply another crucial trick which helps us both at the non-archimedean places and the archimedean place. Experience has shown (e.g.\ \cite{Mo, BLM, AK}) that reciprocity formulae become much easier if a kind of root number is included. Therefore we first apply the functional equation to one of the factors which introduces the root number artificially, and then start the game of applying summation formulae, which now becomes much easier to deal with. As one ${\rm GL}(2)$ functional equation was already applied at the beginning, the final ${\rm GL}(3)$ Voronoi summation formula has to be reduced to a  ${\rm GL}(1)$ Poisson formula. 
 
\subsection{Statement of the formula} We now proceed to the precise statement of the new reciprocity formula which requires a bit of notation. We specify a class of test functions. Fix some very large constant $A$. A function $h = h_T: \Bbb{C} \rightarrow \Bbb{C}$ depending on a parameter $T > 1$ is called \emph{$T$-admissible} if 
\begin{equation}\label{admi}
\begin{split}
&\text{\textbullet \hspace{0.2cm}  it is even and holomorphic in $|\Im t| \leq A+1$;}\\
& \text{\textbullet \hspace{0.2cm}  $h(i(n-1/2)) = 0$ for all $n \in \Bbb{Z}$ with $|n| \leq A+1$;}\\
& \text{\textbullet \hspace{0.2cm}    for all $j \in \Bbb{N}_0 $ we have $ h^{(j)}(t) \ll_j  e^{-(t/T)^2}
 \left(1 + \frac{|t|^2}{T^2}\right)^{j}\left(\frac{1+|t|}{T}\right)^A (1 + |t|)^{-j}  $.} 
\end{split}
\end{equation}
For the rest of the paper let $q$ be a prime and $\chi$ a primitive character modulo $q$. Let $\mathcal{B}(q, \chi)$ denote an orthogonal basis of cuspidal Maa{\ss} newforms (since $\chi$ is primitive, there are no oldforms). For $f \in \mathcal{B}(q, \chi)$ we denote by $t_f \in i\Bbb{R} \cup [-i\vartheta_0, i\vartheta_0]$ its spectral parameter, by $\lambda_f(n)$ its $n$-th Hecke eigenvalue and by $\epsilon_f$ its eigenvalue under the involution $f \mapsto f(-\bar{z})$. This involution flips positive and negative Fourier coefficients \cite[(4.70)]{DFI}.  We denote the two (families of) Eisenstein series corresponding to the cusps $\mathfrak{a} = 0, \infty$ by $\mathscr{E}_{\mathfrak{a}}(t)$ for $t\in \Bbb{R}$. They have Hecke eigenvalues  (\cite[(6.17)]{DFI})
\begin{equation}\label{lambda-eis}
\lambda_{\infty, t}(n) = \sum_{ad = n} \chi(a) (a/d)^{it}, \quad  \lambda_{0, t}(n) = \sum_{ad = n} \chi(d) (a/d)^{it},  
\end{equation}
 with corresponding $L$-functions  
 $$L(s, \mathscr{E}_{\mathfrak{a}}(t)) := \sum_{n} \frac{\lambda_{\mathfrak{a}, t}(n)}{n^s} = \begin{cases} L(s - it, \chi) \zeta(s+it), & \mathfrak{a} = \infty,\\
 L(s + it, \chi) \zeta(s-it), & \mathfrak{a} = 0.\end{cases}$$
 Let $\ell \in \Bbb{N}$ with $(\ell, q) = 1$. For $w\in \Bbb{C}$,   $f \in \mathcal{B}(q, \chi)$ and $\mathfrak{a} \in \{\infty, 0\}$ define
 $$\Lambda_f(\ell; w) = \sum_{\ell_1\ell_2 = \ell} \frac{\mu(\ell_1)\chi(\ell_1) \lambda_f(\ell_2)}{\ell_1^w}, \quad \Lambda_{\mathfrak{a}, t}(\ell; w) = \sum_{\ell_1\ell_2 = \ell} \frac{\mu(\ell_1)\chi(\ell_1) \lambda_{\mathfrak{a}, t}(\ell_2)}{\ell_1^w}.$$
Finally let $F$ be an automorphic form for ${\rm SL}(3, \Bbb{Z})$ with Fourier coefficients 
\begin{equation}\label{eisen-hecke} 
  A_F(n_1, n_2) = A(n_1, n_2) =  \sum_{d \mid (n_1, n_2)} \mu(d) A(n_1/d, 1) A(1, n_2/d),
\end{equation}   
   which can be cuspidal or the minimal parabolic Eisenstein series   ${\tt E}_{0}$   with Hecke eigenvalues
\begin{equation}\label{etau}
A(n, 1)  =\tau_3(n) =  \sum_{abc = n}1. 
\end{equation}
 We denote the archimedean Langlands parameters of $F$ by $\mu = (\mu_1, \mu_2, \mu_3)$ satisfying  $\mu_1 + \mu_2 + \mu_3 = 0$ and write $\theta = \theta_F = \max_j |\Re \mu_j| \leq 5/14$, an admissible exponent for the Ramanujan conjecture on ${\rm GL}(3)$ (not to be confused with $\vartheta \leq 7/64$).  For $F = {\tt E}_{0}$ we have $\theta = 0$.   We write $\tilde{F}$ for the dual form with $A_{\tilde{F}}(n_1, n_2) = A_{F}(n_2, n_1) = \overline{A_F(n_1, n_2)}$. We fix $F$ once and for all, and all implied constants may depend on $F$.

For $s, w \in \Bbb{C}$  and   a $T$-admissible function $h$ we define
\begin{equation}\label{curly-N}
\begin{split}
\mathcal{N}^{\text{cusp}}_{q, \ell}(s, w; h) &:= \sum_{f \in \mathcal{B}(q, \chi)} \epsilon_f \lambda_f(q) \frac{L(s, f \times F) L^{(q)}(w, \bar{f})}{L(1, \text{Ad}^2f)} \frac{\overline{\Lambda_f}(\ell; w)}{\ell^w} h(t_f),\\
\mathcal{N}^{\text{Eis}}_{q, \ell}(s, w; h)& :=\sum_{\mathfrak{a} \in \{\infty, 0\}} \int_{\Bbb{R}}  \frac{L(s,\mathscr{E}_{\mathfrak{a}}(t) \times F) L^{(q)}(w,\mathscr{E}_{\mathfrak{a}}(t))}{L(1 - 2 it, \chi_{\mathfrak{a}})L(1 + 2 it, \overline{\chi_{\mathfrak{a}}})} \frac{\overline{\Lambda_{\mathfrak{a}, t}}(\ell; w)}{\ell^w} h(t) \frac{dt}{2\pi},\\
\end{split}
\end{equation}
where here and henceforth a superscript $^{(N)}$ denotes removal of the Euler factor at primes dividing $N$, $\chi_{\mathfrak{a}} = \chi$ if $\mathfrak{a} = \infty$ and $\chi_{\mathfrak{a}} = \bar{\chi}$ if $\mathfrak{a} = 0$ and $\Re s \not= 1 \not= \Re w$ in the second expression if $F = {\tt E}_{0}$ is not cuspidal. 
 In this case we have  
 \begin{equation}\label{eis-expl}
 L(s,\mathscr{E}_{\mathfrak{a}}(t) \times {\tt E}_0) L^{(q)}(w,\mathscr{E}_{\mathfrak{a}}(t)) =   \big(L(s \mp it , \chi) \zeta(s \pm it )\big)^3  L(w \mp it, \chi) \zeta^{(q)}(w\pm it) 
 \end{equation}
 with the upper sign for $\mathfrak{a} = \infty$ and the lower sign for $\mathfrak{a} = 0$. 
It follows from Lemma \ref{ana-cont}a that $\mathcal{N}^{\text{Eis}}_{q, \ell}(s, w; h)$ as a function of $(s, w)$ with $\Re s, \Re w > 1$  can be continued analytically to a large tube, and its analytic continuation in $\Re s, \Re w < 1$ equals $\mathcal{N}^{\text{Eis}}_{q, \ell}(s, w; h)$ plus some polar terms. 
We write
$$ \mathcal{N}_{q, \ell}(s, w; h)   =  \mathcal{N}^{\text{cusp}}_{q, \ell}(s, w; h) +\mathcal{N}^{\text{Eis}}_{q, \ell}(s, w; h).$$
Compared to \eqref{proto}, the role of $\lambda_f(\ell)$ is played by $\Lambda_f(\ell;w )$ (a convolution of $\lambda_f(\ell)$ which can be recovered by M\"obius inversion), while $\epsilon_f \lambda_f(q)$  resembles the root number from an artificial application of the functional equation. 

For the other side of the formula, we need slightly different mean values. For $T \geq 1$ we call a pair of functions $\mathfrak{h} = (h, h^{\text{hol}}) : (\Bbb{R} \cup [-i\vartheta_0, i\vartheta_0]) \times 2\Bbb{N} \rightarrow \Bbb{C}$ depending on $T$ \emph{weakly $T$-admissible} if it satisfies
\begin{equation}\label{weakly}
  h(t) \ll T^{1 + 2\vartheta_0+\varepsilon}(1 + |t|)^{-20}, \quad t \in \Bbb{R} \cup [-i\vartheta_0, i\vartheta_0], \quad\quad h^{\text{hol}}(k) \ll Tk^{-20}, \quad k \in 2\Bbb{N}
\end{equation}
for all $\varepsilon > 0$.  

Let $\mathcal{B}^{\ast}(N, \text{triv})$ resp.\ $\mathcal{B}_{\text{hol}}^{\ast}(N, \text{triv})$ denote an orthonormal basis of Maa{\ss} resp.\  holomorphic \emph{new}forms of level $N \in \Bbb{N}$ and trivial central character. As before we denote the spectral parameter of $f \in\mathcal{B}^{\ast}(N, \text{triv})$ by $t_f$ and the weight of  $f \in\mathcal{B}_{\text{hol}}^{\ast}(N, \text{triv})$ by $k_f$. Let $\mathfrak{h} = (h, h^{\text{hol}})$ be weakly $T$-admissible for some $T$. 
For $s, w$ with $\Re s, \Re w > 0$ and  $N \in \Bbb{N}$ with $q^2 \nmid N$ define
\begin{equation}\label{curly-M}
\begin{split}
\mathcal{M}^{\text{Maa{\ss}}, \pm}_{N}(s, w; \mathfrak{h}) &:= \sum_{N_0 \mid N} \sum_{f \in \mathcal{B}^{\ast}(N_0, \text{triv})} \epsilon_f^{(1\mp 1)/2}   \frac{L(s, f \times \tilde{F}) L(w, f\times \bar{\chi})}{L(1, \text{Ad}^2f)}   \tilde{L}_N(s, w; f)h(t_f),\\
\mathcal{M}^{\text{hol} }_{N}(s, w; \mathfrak{h}) &:= \sum_{N_0 \mid N} \sum_{f \in \mathcal{B}_{\text{hol}}^{\ast}(N_0, \text{triv})}   \frac{L(s, f \times \tilde{F}) L(w, f\times \bar{\chi})}{L(1, \text{Ad}^2f)}  \tilde{L}_N(s, w; f) h(k_f),\\
\end{split}
\end{equation}
\begin{displaymath}
\begin{split}
\mathcal{M}^{\text{Eis}}_{N}(s, w; \mathfrak{h})& :=\sum_{\psi : c_{\psi}^2 \mid N}  \int_{\Bbb{R}}   \frac{L(s + it, \tilde{F} \times \psi) L(s - it, \tilde{F} \times \bar{\psi}) L(w + it, \bar{\chi}\psi)L(w - it, \bar{\chi}\bar{\psi})}{L(1 + 2 it, \psi^2)L(1 - 2it, \bar{\psi}^2)} \\
&\quad\quad\quad\quad\quad\tilde{L}_N(s, w; (t, \psi)) h(t) \frac{dt}{2\pi},\\
\end{split}
\end{displaymath}
where the $\psi$-sum runs over all primitive Dirichlet characters $\psi$ whose conductor $c_{\psi}$ satisfies $c_{\psi}^2 \mid N$ (in particular, only the trivial character contributes if $N$ is squarefree) and we assume $\Re s \not= 1 \not= \Re w$ in the expression $\mathcal{M}^{\text{Eis}}_{N}(s, w; h)$ if $F = {\tt E}_{0}$. The expressions $\tilde{L}_N(s, w; f)$ and $\tilde{L}_N(s, w; (t, \psi))$ are  Euler polynomials (local correction factors) that are defined in \eqref{cusp-corr} and \eqref{Eis-corr}. We only need to know
\begin{equation}\label{corr-bound-final}
\tilde{L}_{N}(s, w; f), \, \tilde{L}_N(s, w; (t, \psi)) \ll N^{\theta - 1}(qN)^{\varepsilon}
\end{equation}
uniformly in $\Re s, \Re w\geq 1/2$, and for all $f \in \mathcal{B}^{\ast}(N_0, \text{triv})$ for $N_0 \mid N$ and all  $t\in \Bbb{R}$, $q^2 \nmid N$ and primitive Dirichlet characters $\psi$ of conductor $c_{\psi}$ with $c_{\psi}^2 \mid N$. 

 For $F = {\tt E}_{0}$, the term $\mathcal{M}^{\text{Eis}}_{N}(s, w; \mathfrak{h})$ for $\Re s, \Re w > 1$ can be continued to $\Re s, \Re w\geq 1/2$ (under suitable conditions on $h$), and its analytic continuation in $\Re s, \Re w < 1$ equals $\mathcal{M}^{\text{Eis}}_{N}(s, w; \mathfrak{h})$ plus some polar terms, cf.\ Lemma \ref{ana-cont}b.  We write
$$\mathcal{M}^{\pm}_{N}(s, w; \mathfrak{h}) = \mathcal{M}^{\text{Maa{\ss}}, \pm}_{N}(s, w; \mathfrak{h}) + \delta_{\pm = +}\mathcal{M}^{\text{hol} }_{N}(s, w; \mathfrak{h}) + \mathcal{M}^{\text{Eis}}_{N}(s, w; \mathfrak{h}).$$
We are now ready to state the reciprocity formula.

\begin{theorem}\label{thm2} Let $h$ be a $T$-admissible function for some $T \geq 1$. Let $q$ be prime, $\chi$ a primitive character modulo $q$ with Gau{\ss} sum $\tau(\chi)$  and $(\ell, q) = 1$. Let $s, w \in \Bbb{C}$ be such that $1/2 \leq \Re s \leq \Re w < 1$, and define $s' = (s+w)/2$, $w' = 1+w/2 - 3s/2$. Then
\begin{displaymath}
\begin{split}
  \mathcal{N}_{q, \ell}(s, w; h) = \mathcal{G}_{q, \ell}(s, w; h) + \frac{   \tau(\chi)^2}{q^{(1+ 3s'-3w')/2}}  \sum_{\pm} \sum_{\epsilon \in \{\pm 1\}} \left(\mathcal{M}^{\pm \epsilon}_{\ell}(s', w';\mathcal{T}_{s', w'}^{\pm, \epsilon}h) - \mathcal{M}^{\pm \epsilon}_{\ell q}(s', w'; \mathcal{T}_{s', w'}^{\pm, \epsilon}h) \right).
 \end{split}
 \end{displaymath} 
The transform $\mathcal{T}_{s', w'}^{\pm, \epsilon}$ is defined in \eqref{defHtrafo} and \eqref{defT} and is weakly $T$-admissible in the sense of \eqref{weakly} for $\Re s= \Re w = 1/2$.  The main term  $\mathcal{G}_{q, \ell}(s, w; h)$ is defined in \eqref{defGmain} and vanishes unless $F = {\tt E}_0$.  
 If $F = {\tt E}_{0}$   and $h$ has   simple     zeros at $\pm i(w-1)$ and   triple zeros at  $\pm i(s-1)$, then 
 \begin{equation}\label{bound-G}
 \mathcal{G}_{q, \ell}(s, w; h) \ll_{\varepsilon}  T^2 q \ell^{-1}  (q\ell T)^{\varepsilon}
 \end{equation}
 for  $ \Re s = \Re w =1/2$ and any $\varepsilon > 0$.  \end{theorem}

The additional zeros of $h$ required in  \eqref{bound-G} play the same role as in \cite{BHM}: we are ultimately interested in bounding a spectral sum over cusp forms (cf.\ Theorem \ref{thm3}), but the artificially added Eisenstein term produces a much bigger term, coming from additional polar terms. The zeros of $h$ make the corresponding residues vanish. 

A quantitative version of Theorem \ref{thm2} (without root numbers)  reads as follows. For $\tau \in \Bbb{R}$, $t \in \Bbb{R} \cup [-i\vartheta_0, i\vartheta_0]$ define
$$Q_{\tau}(t) =  \frac{ 1 -   i \sinh(\pi \tau)/\cosh(\pi t)}{1 + i \sinh(\pi \tau)/\cosh(\pi t) } = 1 + O_{\tau}(e^{-\pi |t|}).$$

\begin{theorem}\label{thm3} Let $q$ be a prime and $\chi$ a primitive character modulo $q$. Let $\ell \in \Bbb{N}$ be coprime to $q$, and let $T\geq 1$, $\tau\in \Bbb{R}$ and $h$ a $T$-admissible function with the additional property that it has a   zero at $\pm i(1/2 \pm i\tau \pm n)$ for all sign combinations and all $n \in \Bbb{N}$, $1 \leq n \leq A+1$, and a quadruple zero at $\pm i(1/2 \pm i\tau)$ for all sign combinations.  
 Then
\begin{displaymath}
\begin{split}
\sum_{ f \in \mathcal{B}(q, \chi) }&\begin{cases} 1, & \epsilon_f = 1\\ Q_{\tau}(t_f), & \epsilon_f = -1\end{cases} \Bigg\} \lambda_f(\ell) \frac{|L(1/2+i\tau , f)|^4 }{L(1, \text{{\rm Ad}}^2 f)} h(t_f)\\
& \ll_{\tau, \varepsilon} \left( qT^2 \ell^{-1/2} + T^2 q^{1/2}+ q^{3/4} T^{7/4}+ T^{1+2\vartheta_0}q^{1/2} (q^{3/8} \ell^{3/4} + q^{1/4} \ell  )\right)(Tq\ell)^{\varepsilon}
\end{split}
\end{displaymath}
for all $\varepsilon > 0$. 
\end{theorem}
A typical candidate for $h$ is
\begin{equation}\label{defh}
\begin{split}
h(t) = &e^{-(t/T)^2}T^{-6[A]-24}    (t^2 + (\textstyle\frac{1}{2} +   i \tau)^2)^3(t^2 + (\textstyle\frac{1}{2} -   i \tau)^2)^3 \\&\displaystyle\prod_{n=1}^{[A]+2}(t^2 + (n -  \textstyle\frac{1}{2})^2)(t^2 + (n +i\tau -  \textstyle\frac{1}{2})^2)(t^2 + (n - i\tau -  \textstyle\frac{1}{2})^2),
\end{split}
\end{equation}
a function that is positive on $\Bbb{R} \cup [-i\vartheta_0, i\vartheta_0]$ and satisfies $h(t) \asymp 1$ for $t \asymp T$. 
Theorem \ref{thm1} can be obtained easily from Theorem \ref{thm3} by standard amplification. 

\subsection{Heuristics}\label{heu} For the reader's convenience we give a short heuristic overview on the genesis of the reciprocity formula \eqref{proto} and its companion in Theorem \ref{thm2} in the case of $F = {\tt E}_0$. We have
\begin{displaymath}
\begin{split}
  \sum_{f \in \mathcal{B}(q, \chi)} |L(1/2, f)|^4 \lambda_f(\ell) \approx \sum_f  \sum_{n_1, n_2, m_1, m_2 \asymp q^{1/2}} \frac{ \lambda_f(n_1n_2\ell) \overline{\lambda}_f(m_1m_2)}{(n_1n_2m_1m_2)^{1/2}}
  \end{split}
\end{displaymath}
We apply the Kuznetsov formula, getting an off-diagonal term roughly of the form
$$  \sum_{q \mid c \asymp q\ell^{1/2}} \frac{1}{c}  \sum_{n_1, n_2, m_1, m_2 \asymp q^{1/2}}    S_{\chi}(n_1n_2 \ell, m_1m_2, c).$$
We apply Poisson summation in $n_2, m_1, m_2$ getting roughly
$$\frac{1}{q^{1/2}\ell} \sum_{q \mid c \asymp q\ell^{1/2}} \sum_{n_1 \asymp q^{1/2}} \sum_{n_2, m_1, m_2 \asymp (q\ell)^{1/2}} \chi(\overline{n_1\ell} n_2) e\left(\frac{\overline{n_1\ell} n_2m_1m_2}{c}\right).$$
We apply the additive reciprocity formula, replacing  $e( \overline{n_1\ell} n_2m_1m_2/c)$ with $e( \bar{c}n_2m_1m_2/n_1\ell)$, and apply  Poisson in $n_2, m_1, m_2$ again, getting
$$\frac{1}{q^{1/2}\ell} \sum_{q \mid c \asymp q\ell^{1/2}}  \sum_{n_1 \asymp q^{1/2}} \sum_{n_2 \asymp q\ell^{1/2}} \sum_{m_1, m_2 \asymp \ell^{1/2}} \bar{\chi}(n_2) S(cm_1m_2, n_2\bar{q}, n_1\ell).$$
We write $c = q\gamma$ and cancel $q$ in the Kloosterman sum. Then we apply the Kuznetsov formula backwards and recognize central values of $L$-values, namely
\begin{equation}\label{arrive}
\frac{q^{1/2}}{\ell^{1/2}} \sum_{f \in \mathcal{B}(\ell, \text{triv})} L(1/2, f)^3 L(1/2, f \times \bar{\chi}).
\end{equation}

Alternatively, we can apply the functional equation to one of the $L$-factors, so that we are looking at something roughly of the shape
$$ \sum_{f \in \mathcal{B}(q, \chi)} \lambda_f(q) L(1/2, f)^3 L(1/2, \bar{f}) \overline{\lambda_f}(\ell)$$
where for convenience we replaced $\lambda_f(\ell)$ with $ \overline{\lambda_f}(\ell)$. Now the off-diagonal term after Kuznetsov is roughly
\begin{displaymath}
\begin{split}
& \sum_{q \mid c \asymp q^{3/2} \ell^{1/2}} \sum_{n_1, n_2, n_3, m \asymp q^{1/2}} \frac{S_{\chi}(m\ell, q n_1 n_2 n_3, c)}{c} \approx \frac{1}{q^{1/2}}\sum_{\gamma \asymp (q \ell)^{1/2}} \sum_{n_1, n_2, n_3, m \asymp q^{1/2}} \frac{S (m\ell\bar{q},  n_1 n_2 n_3, \gamma)\chi(\gamma \overline{m})}{\gamma}
\end{split}
\end{displaymath}
by twisted multiplicativity. We apply Poisson in $n_1, n_2, n_3$ getting roughly
$$\frac{1}{\ell}\sum_{\gamma \asymp (q \ell)^{1/2}} \sum_{ m \asymp q^{1/2}}\sum_{n_1, n_2, n_3 \asymp \ell^{1/2}}  \chi(\gamma \overline{m }) e\left( \frac{\overline{m\ell} n_1n_2n_3 q}{\gamma}\right).$$
We apply reciprocity and   Poisson in $\gamma$ getting
$$\frac{1}{\ell  q^{1/2}}\sum_{\gamma \asymp q \ell^{1/2}}  \sum_{ m \asymp q^{1/2}}\sum_{n_1, n_2, n_3 \asymp \ell^{1/2}}  \bar{\chi}(\gamma) S(n_1n_2n_3, \gamma, \ell m),$$
which leads to \eqref{arrive} after an application of the Kuznetsov formula. 

Both approaches are essentially equivalent, but the second one is technically easier, which is why we take this route. 

\subsection{Notation} As usual, for $q \mid c$ we write
$$S_{\chi}(m, n, c) = \underset{d\, (\text{mod } c)}{\left.{\sum}\right.^{\ast}} \chi(d) e\left(\frac{md + \bar{d}n}{c}\right)$$
for the twisted Kloosterman sum, and we write 
$\tau(\chi) = S_{\chi}(1, 0, q)$ 
for the standard Gau{\ss} sum. By $v_p(n)$ we denote the $p$-adic valuation of $n$. The statement $a\mid b^{\infty}$ means that $a$ has only  prime divisors that divide $b$. For $n \in \Bbb{N}$ we denote by $\omega(n)$ the number of distinct prime divisors of $n$. By an $\varepsilon$-neighbourhood of a strip $c_1  \leq \Re s \leq c_2$ (possibly $c_2 = \infty$) we mean the open set $c_1 - \varepsilon <\Re s< c_2+\varepsilon$ for some sufficiently small $\varepsilon > 0$, and similarly for multidimensional tubes. It is convenient to introduce the function
\begin{equation}\label{cs}
c(s) = \max\big(0, 1/2 - \Re s, (1-\Re s)/2\big)
\end{equation}
so that the convexity bound for Dirichlet $L$-functions states 
$L(s, \chi) \ll_s q^{c(s) + \varepsilon}.$
The symbol $\varepsilon$ denotes an arbitrarily small positive real number whose value may change from line to line. 


 \section{Versions of the Kuznetsov formula}
As our reciprocity formula is not symmetric, we need two different versions of the Kuznetsov formula, one with central character and another one without central character. We start with the latter and quote from \cite[Section 3]{BK}. 
 
Let $N \in \Bbb{N}$ and
$$ N \nu(N) := [\Gamma_0(1) : \Gamma_0(N)], \quad \text{i.e.} \quad \nu(N) = \prod_{p \mid N} \left(1 + p^{-1}\right).$$ Eisenstein series for $\Gamma_0(N)$ are parametrized by a continuous parameter $s = 1/2 + it$, $t \in \Bbb{R}$, and pairs $(\psi, M)$ where $\psi$ is a primitive Dirichlet character of conductor $c_{\psi}$ and $M \in \Bbb{N}$ satisfies $c_{\psi}^2 \mid M \mid N$. We define  
\begin{equation*}
\mathfrak{n}^2(M) :=  \frac{1}{M} \prod_{\substack{p \mid N\\ p \nmid (M, N/M)}} \frac{p}{(p+1)} \prod_{ p \mid (M, N/M)} \frac{p-1}{p+1} =: \frac{1}{M} \tilde{ \mathfrak{n}}^2(M)
\end{equation*}
and write
$$ M = c_{\psi} M_1 M_2, \quad \text{where} \quad (M_2, c_{\psi}) = 1, \quad M_1 \mid c_{\psi}^{\infty}. $$
The normalized Eisenstein series $E_{\psi, M, N}(z, s)$ of level $N$ corresponding to   $(\psi, M)$ has the Fourier expansion $$E_{\psi, M, N}(z, 1/2 + it) = \rho^{(0)}_{\psi, M, N}(t, y) + \frac{2 \pi^{1/2 + it} y^{1/2}}{\Gamma(1/2 + it)} \sum_{n \not= 0} \rho_{\psi, M, N}(n, t) K_{it}(2 \pi |n| y) e(nx),$$ where for $n \not= 0$ we have 
  \begin{equation}\label{rho-eis}
    \rho_{\psi, M, N}(n, t)  = \frac{C(\psi, M, t)|n|^{it}}{(N \nu(N))^{1/2}\tilde{ \mathfrak{n}}(M) L^{(N)}(1 + 2it, \psi^2)}  \left(\frac{M_1}{M_2}\right)^{1/2}  \sum_{\delta \mid M_2} \delta \mu(M_2/\delta)  \bar{\psi}(\delta)  \sum_{\substack{cM_1\delta f =  n \\ (c, N/M) = 1}}\frac{\psi(c)}{c^{2it}}   \bar{\psi}(f)
    \end{equation}
    for a constant $C(\psi, M, t)$ with $|C(\psi, M, t)| = 1$. 
    
    The cuspidal spectrum is parametrized by  pairs $(f, M)$ of $\Gamma_0(N)$-normalized newforms $f$ of level $N_0 \mid N$ and integers $ M \mid N/N_0$. If $f$ is a Maa{\ss} form with spectral parameter $t$ and Hecke eigenvalues $\lambda_f(n)$, we write the Fourier expansion of the pair $(f, M)$ as
$$ (2 \cosh(\pi t) y)^{1/2} \sum_{n \not= 0} \rho_{f, M, N}(n)   K_{it}(2 \pi |n| y) e(nx)$$
 with 
 \begin{equation}\label{rho-cusp}
 \rho_{f, M, N}(n) = \frac{1}{L^{\ast}(1, \text{Ad}^2 f)^{1/2}(N\nu(N))^{1/2}} \prod_{p \mid N_0} \left(1 - \frac{1}{p^2}\right)^{1/2}  \sum_{d \mid M} \xi_{f}(M, d) \frac{d}{M^{1/2}} \lambda_f(n/d)  
 \end{equation}
 for $n \in \Bbb{N}$ with the convention $\lambda_f(n) = 0$ for $n \not \in \Bbb{Z}$. Here $\xi_f(M, d)$ is a certain arithmetic function defined in \cite[Section 3.2]{BK}, and all we need to know is the bound \cite[(3.9)]{BK}
\begin{equation}\label{xi-arithmetic}
\xi_{f}(M, d)  \ll M^{\varepsilon}(M/d)^{\vartheta}. 
\end{equation}
Moreover, we use the   notation 
$$L^{\ast}(s, \text{Ad}^2 f) = \zeta^{(N)}(2s) \sum_n \lambda_f(n^2)n^{-s} $$
(which may differ from $L(s, \text{Ad}^2f)$ by finitely many Euler factors if $N$ is not squarefree, cf.\ \cite{Li}).  For $-n \in \Bbb{N}$ we have $\rho_{f, M, N}(n) = \epsilon_f \rho_{f, M, N}(-n)$ where $\epsilon_f \in \{\pm 1\}$ is the parity of $f$.
 
 If $f$ is a holomorphic newform of weight $k$ and level $N_0 \mid N$, we write the Fourier expansion of the pair $(f, M)$ as
$$\left(\frac{2\pi^2}{\Gamma(k)}\right)^{1/2}\sum_{n > 0} \rho_{f, M, N}(n) (4\pi n)^{(k-1)/2} e(nz).$$
Then \eqref{rho-cusp} remains true for  $n \in \Bbb{N}$, and  for negative $n$ we define $\rho_{f, M, N}(n)  = 0$. 

For $x > 0$ we define the integral kernels 
\begin{equation}\label{kernel}
\begin{split}
& \mathcal{J}^+(x, t) := \frac{\pi  i}{ \sinh(\pi t)} (J_{2 it}(4 \pi x) - J_{-2it}(4 \pi x)), \\
& \mathcal{J}^-(x, t) := \frac{\pi  i}{ \sinh(\pi t)} (I_{2 it}(4 \pi x) - I_{-2it}(4 \pi x)) = 4 \cosh(\pi t) K_{2it}(4\pi x),\\
& \mathcal{J}^{\text{hol}}(x, k) := 2\pi i^k J_{k-1}(4\pi x) =  \mathcal{J}^+(x, (k-1)/(2i)) , \quad k \in 2\Bbb{N}. 
 \end{split}
 \end{equation} 
If $H \in C^3((0, \infty))$ satisfies  $x^j H^{(j)}(x) \ll \min(x, x^{-3/2})$ for $0 \leq j \leq 3$ and $n, m, N \in \Bbb{N}$, then  we have 
 \begin{equation}\label{kuz2} 
\begin{split}
\sum_{N \mid c} &\frac{S(\pm n, m, c)}{c}H \left(\frac{\sqrt{nm}}{c}\right) =  \mathcal{A}_N(\pm n, m; \mathscr{L}_{\pm} H)   \end{split}
\end{equation}
where 
\begin{equation}\label{Hback}
\begin{split}
\mathscr{L}_{\pm}H = (\mathscr{L}^{\pm}H, \mathscr{L}^{\text{hol}}H), \quad \mathscr{L}^{\diamondsuit}H = \int_{0}^{\infty} \mathcal{J}^{\diamondsuit}(x, .) H(x) \frac{dx}{x}. 
\end{split}
\end{equation}
for $\diamondsuit \in \{+, -, \text{hol}\}$ and 
$$\mathcal{A}_N(n, m; \mathfrak{h})  :=   \mathcal{A}_N^{\text{Maa{\ss}}}(n, m; \mathfrak{h}) +  \mathcal{A}_N^{\text{Eis}}(n, m; \mathfrak{h})+   \mathcal{A}_N^{\text{hol}}(n, m; \mathfrak{h})$$
 for $n, m \in \Bbb{Z} \setminus \{0\}$ with
  \begin{equation}\label{defA}
\begin{split}
& \mathcal{A}_N^{\text{Maa{\ss}}}(n, m; \mathfrak{h}) := \sum_{N_0M \mid N}   \sum_{f \in \mathcal{B}^{\ast}(N_0, \text{triv})} \rho_{f, M, N}(n) \overline{ \rho_{f, M, q}(m)} h(t_f),\\
&  \mathcal{A}_N^{\text{Eis}}(n, m; \mathfrak{h}) := \sum_{c_{\chi}^2\mid M \mid N} \int_{\Bbb{R}} \rho_{\chi, M, N}(n, t) \overline{ \rho_{\chi, M, q}(m, t)} h(t) \frac{dt}{2\pi},\\
& 
\mathcal{A}_N^{\text{hol}}(n, m; \mathfrak{h}) :=  \sum_{N_0M \mid N}   \sum_{ f \in \mathcal{B}_{\text{hol}}^{\ast}(N_0, \text{triv}) } \rho_{f, M, N}(n) \overline{ \rho_{f, M, q}(m)} h^{\text{hol}}(k_f)
 \end{split}
\end{equation}
with with $\mathfrak{h} = (h, h^{\text{hol}})$ and the notation $\mathcal{B}^{\ast}(N_0, \text{triv})$ and $\mathcal{B}_{\text{hol}}^{\ast}(N_0, \text{triv})$ as in \eqref{curly-M}. 

We need the Kuznetsov formula for prime level $q$ and primitive central character $\chi$ in the other direction. The notation is a little easier here, because there are no cuspidal oldforms and also the Eisenstein spectrum contains in some sense only newforms, in particular the classical parametrization in terms of the cusps $\infty$ and $0$ and the adelic parametrization in terms of the two pairs $(\chi, \text{triv})$, $( \text{triv}, \chi)$ is the same. Following \cite[Section 5]{KL} (or \cite[Section 7]{DFI}), the Fourier coefficients $\rho_{\mathfrak{a}, t}(n)$ are 
$$\rho_{\mathfrak{a}, t}(n) = \frac{C(\mathfrak{a}, t)\lambda_{\mathfrak{a}, t}(|n|)}{q^{1/2} L(1 + 2 it, \overline{\chi_{\mathfrak{a}}})}$$
for some constant $C(\mathfrak{a}, t)$ with $|C(\mathfrak{a}, t)| = 1$,  $\lambda_{\mathfrak{a}, t}$ as in \eqref{lambda-eis} and $\chi_{\mathfrak{a}} = \chi$ if $\mathfrak{a} = \infty$ and $\chi_{\mathfrak{a}} = \bar{\chi}$ if $\mathfrak{a} = 0$  Similarly, for $f \in \mathcal{B}(q, \chi)$ we use that (cf.\ \cite[Section 2]{Li})
$$\underset{s=1}{\text{res}} \zeta(2s) \sum_n \frac{|\lambda_f(n)|^2}{n^s} = \underset{s=1}{\text{res}}L(s, f \times \bar{f}) (1 - q^{-2s})^{-1} (1 - q^{-s}) = L(1, \text{Ad}^2f) \frac{q}{q+1}$$
to conclude that the Fourier coefficients are
$$\rho_f(n) = \frac{\lambda_f(n)}{q^{1/2} L(1, \text{Ad}^2f)}$$
for $n \in \Bbb{N}$ and $\rho_f(-n) = \epsilon_f \rho_f(n)$. We can now state the Kuznetsov formula directly in terms of Hecke eigenvalues, and we only need the ``opposite sign" case: for $n, m > 0$ and a $T$-admissible function $h$ we have 
\begin{equation}\label{kuz-neben}
\begin{split}
&\frac{1}{q }\sum_{f \in \mathcal{B}(q, \chi)} \epsilon_f  \frac{\lambda_f(n)\overline{\lambda_f(m)}}{L(1, \text{Ad}^2 f)} h(t_f) + \frac{1}{q } \sum_{\mathfrak{a} \in \{\infty, 0\} }\int_{\Bbb{R}}   \frac{\lambda_{\mathfrak{a}, t}(n) \overline{\lambda_{\mathfrak{a}, t}(m)}}{L(1 - 2it, \chi_{\mathfrak{a}})L(1 + 2it, \overline{\chi_{\mathfrak{a}}})}  h(t) \frac{dt}{2\pi} \\
&= \sum_{q \mid c} \frac{S_{\chi}(-m, n, c)}{c} \mathscr{K}h \left(\frac{\sqrt{nm}}{c}\right)
\end{split}
\end{equation}
where
\begin{equation*}
 \mathscr{K}h(x) :=   \int _{-\infty}^{\infty} \mathcal{J}^-(x, t) h(t)   t \tanh(\pi t) \frac{dt}{2\pi^2}.
  \end{equation*}
(Note that $\chi$ is even, so that $S_{\chi}(-m, n, c) = S_{\chi}(m, -n, c)$  holds.)

\section{Analytic lemmas} 

 Our first lemma analyzes the integral transform in the opposite sign Kuznetsov formula \eqref{kuz-neben}. 
 
 \begin{lemma}\label{BLM-lemma} Let $T \geq 1$ and let $h$ be a $T$-admissible function. Let $B > 0$, $j \in \Bbb{N}$. If $A$ in \eqref{admi} is sufficiently large in terms of $B$ and  $j$, then 
\begin{equation}\label{bound-BLM}
x^j \frac{d^j}{dx^j} \mathscr{K}h (x) \ll_{B, j} T \min\left( (x/T)^B, (x/T)^{-B}\right).  
\end{equation}
 The Mellin transform $\widehat{\mathscr{K}h }$ is holomorphic in $|\Re u | < B$ and satisfies
\begin{equation}\label{coro-bound}
\widehat{\mathscr{K}h}(u) \ll T^{1+\Re u} (1 + |u|)^{-B}
\end{equation}
 in this region. 
\end{lemma}

\textbf{Proof}  We follow \cite[Lemma 4]{BLM} and consider three ranges.

\emph{Range I:} $x \leq 1$.  We express the kernel $\mathcal{J}^{-}(x, t)\tanh(\pi t)$ in terms of Bessel-$I$-functions  divided by $\cosh(\pi t)$  (cf.\ \eqref{kernel})
and differentiate under the integral sign using \cite[(A.2)]{BLM}
 \begin{equation*}
    I^{(j)}_{2it}(x)  =   \left(\frac{1}{2}\right)^j \sum_{n=0}^j \left(\begin{array}{c} j\\ n\end{array}\right) I_{2it - j + 2n}(x).
 \end{equation*}
 In those terms containing $I_{2it + m}(x)/\cosh(\pi t)$ for $|m| \leq j$ we shift the contour down to $\Im t = -A$, while in those terms containing $I_{-2it + m}(x)/\cosh(\pi t)$ for $|m| \leq j$ we shift the contour up to $\Im t = A$. The zeros of $h$ cancel the poles of $ \cosh(\pi t)^{-1}$. Estimating trivially using \cite[(A.6)]{BLM}
 \begin{equation*}
   e^{-\pi |t|}  I_{2it}(x)     \ll_{\Im t}   \frac{x^{-2\Im t}}{(1+|t|)^{1/2 - 2\Im t}} \quad (0 < x \leq 1), 
 \end{equation*}
we obtain
 $$x^j \frac{d^j}{dx^j} \mathscr{K}h (x)\ll x^j \int_{\Bbb{R}} e^{-(t/T)^2}\left(\frac{1 + |t|}{T}\right)^A \frac{x^{A-j}}{(1+|t|)^{A - j + 1/2}} |t| dt \ll \left(\frac{x}{T}\right)^{A - j - 3/2}$$
 which implies \eqref{bound-BLM} for $x \leq 1$. 
 
 \emph{Range II:} $1 \leq x \leq T^{13/12}$. Let
 $$h_{\text{spec}}(t)   :=   \frac{1}{2\pi^2}  h_T(t) t \tanh(\pi t)$$
which has  Fourier transform
\begin{equation*}
\check{h}_{\text{spec}}(v)   =  \int_{-\infty}^{\infty} h_{\text{spec}}(t) e^{-itv} dt   =  \frac{T^2}{2\pi^2} \int_{-\infty}^{\infty} h(tT) \tanh(\pi t T) e^{-itTv} t \, dt.
\end{equation*}
We have by \eqref{admi} (and the Leibniz rule) 
\begin{displaymath}
\begin{split}
\frac{d^n}{dt^n} \big(h(tT) t\big)&  \ll_n  e^{-t^2} (1 + |t|)^{2n+1}\left(\frac{1+|t|T}{T}  \right)^{A-n} 
\end{split}
\end{displaymath}
for $n\in \Bbb{N}$ and 
$$ \frac{d^n}{dt^n}  \tanh(\pi t T)    \ll_n   \frac{T^n}{\cosh(\pi t T)^2} \ll_n \left(\frac{1+|t|T}{T}  \right)^{-n}  $$
 for $n \geq 1$, so that again by the Leibniz rule
 $$ \frac{d^n}{dt^n}\big(h(tT) \tanh(\pi t T)   t\big) \ll e^{-t^2} (1 + |t|)^{2n+1}\left(\frac{1+|t|T}{T}  \right)^{A-n} $$
 for all $n \in \Bbb{N}$. Integrating by parts $n \leq A$ times, we thus obtain
 \begin{equation*}
 \check{h}^{(j)}_{\text{spec}}(v)  \ll_{j, n} \frac{T^{2+j}}{(1+T|v|)^{A_1}}
 \end{equation*}
for all $j \in \Bbb{N}_0$, $A_1 \leq A$. Instead of partial integration, we can also shift the contour to $\Im t = - \text{sgn}(v) A$. Again the zeros of $h(tT)$ cancel the poles of $\tanh(\pi t T)$. This gives us the alternative bound
 \begin{equation}\label{check2}
 \check{h}^{(j)}_{\text{spec}}(v)  \ll_{j} T^{2+j} e^{-A|v|}.
 \end{equation}
The previous two bounds are \cite[(6.3), (6.4)]{BLM}, and we can now literally follow the argument given there. By Parseval we obtain 
$$
\mathscr{K}h (x) =x^j \frac{d^j}{dx^j} \int_{-\infty}^{\infty} e^{ix\sinh(v/2)}  \check{h}_{\text{spec}}(v) dv,$$
and the second bound \eqref{check2} is only needed to justify integration by parts in this integral.  Using the Taylor expansion of $\sinh(v/2)$ we obtain as in \cite[(6.9), (6.12)-(6.14)]{BLM} that 
\begin{displaymath}
\begin{split}
 x^j \frac{d^j}{dx^j} 
\mathscr{K}h (x) = &\sum_{\alpha=0}^{3A_2} \sum_{\beta= 0}^{\lfloor \alpha/3\rfloor} c_{\alpha, \beta}  x^{\beta}     \sum_{n=0}^j   \sum_{\gamma = n}^{A_3} d_{n, \gamma} \int_{-\infty}^{\infty} e^{ixv/2}  v^{\alpha+\gamma} \check{h}^{(n)}_{\text{spec}}(v) dv \\
& + O\left(T^{2+(j-A_1-3)/4} + T^{-\frac{7}{6}A_2 + \frac{1}{4}} +   T^{(1-3A_3 +4j)/4}\right)
\end{split}
\end{displaymath}
for any constants $A_1, A_2, A_3$ satisfying $A_1 \geq j+2$, $A_1 \geq 3A_2 + A_3 + 2$ and $A_1 \leq A$, and certain absolute constants $c_{\alpha, \beta} $, $d_{n, \gamma}$. If $A$ is sufficiently large in terms of $j$, we can make the error term $\ll T^{-A/10}$ which is acceptable for  \eqref{bound-BLM} in the present range $1 \leq x \leq T^{13/12}$. The main term is by Fourier inversion a linear combination of terms of the form
\begin{displaymath}
\begin{split}
x^{\beta} \frac{d^{\alpha+\gamma}}{dx^{\alpha+\gamma}} \left(x^{n} h_{\text{spec}}(x/2)\right)   \ll   x^{\beta}  e^{-x^2/4T^2} x^n \left( \frac{1 + |x|}{T}\right)^A \left(1 + \frac{|x|}{T}\right)^{2(\alpha + \gamma)}   \frac{1}{|x|^{\alpha+\gamma}}
\end{split}
\end{displaymath}
 with $\alpha \leq 3A_2$, $\beta \leq \alpha/3$, $n \leq j $, $n \leq \gamma \leq A_3$.  If $A$ is sufficiently large, this is easily seen to be acceptable for \eqref{bound-BLM} in the present range $1 \leq x \leq T^{13/12}$. 
 
\emph{Range III:} $x \geq T^{13/12}$. Here the desired bound follows simply from the exponential decay of the Bessel $K$-function. 

This completes the proof of \eqref{bound-BLM}. In order to derive \eqref{coro-bound}, we integrate by parts $[B]+1$ times and estimate trivially using \eqref{bound-BLM} for $|\Re u | < B$ and $|\Im u| \geq 1$, and we estimate trivially without integration by parts for $|\Re u| < B$ and $|\Im u| < 1$. \\

The next lemma analyzes the inverse transform in the Kuznetsov formula \eqref{kuz2}.  We recall the notation \eqref{Hback}. 

\begin{lemma}\label{lem3} 
Let $\rho  > 0$, $B \geq 5$, $c\in \Bbb{R}$, $T \geq 1$, $\varepsilon > 0$ sufficiently small.  Let $\Phi$ be a function that is meromorphic in $-2\rho < \Re u < B$ with at most finitely many poles $u_1, \ldots, u_n$ (listed according to multiplicity), all of which are located in $-2\rho < \Re u < 0$. Suppose that $\Phi$ satisfies $$(u-u_1) \cdot\ldots\cdot  (u-u_n)\Phi(u) \ll T^{c-\Re u} (1 + |u|)^{-B}.$$  We write 
$$\widecheck{\Phi}(x) = \int_{(0)} \Phi(u) x^{-u} \frac{du}{2\pi i}$$
for the inverse Mellin transform. Then $\mathscr{L}^{\pm} \widecheck{\Phi}$ has meromorphic continuation  to $|\Im t | < \rho$ with poles at most at $t = \pm i u_j/2$, $j = 1, \ldots, n$, and it 
satisfies
\begin{equation}\label{L1}
\mathscr{L}^{\pm} \widecheck{\Phi}(t)  \ll_{B, \varepsilon } T^c  (1+|t|)^{-B}, \quad t \in \Bbb{R},
\end{equation}
and
\begin{equation}\label{L2}
(t-\textstyle \frac{1}{2} iu_1)  \cdot\ldots\cdot(t- \textstyle \frac{1}{2} iu_n) \mathscr{L}^{\pm} \widecheck{\Phi}(t) 
 \ll_{B, u_1, \ldots, u_n, t, \varepsilon} 
T^{c+ 2|\Im t|+\varepsilon}, \quad |\Im t| < \rho.
\end{equation}
Moreover
\begin{equation}\label{L3}
 \mathscr{L}^{\text{{\rm hol}}} \widecheck{\Phi}(k) \ll_B T^c k^{-B}, \quad  k \in 2\Bbb{N}.
 \end{equation}
All implied constants are independent of $T$. 
\end{lemma}

\textbf{Proof.} This is almost verbatim   \cite[Lemma 3]{BK}. Using the definition \eqref{Hback}  and exchanging integrals, we have for $t \in \Bbb{R}$ and $\varepsilon > 0$ sufficiently small that 
\begin{equation}\label{start}
\mathscr{L}^{\pm}\widecheck{\Phi}(t) = \int_{(  \varepsilon)}\widehat{ \mathcal{J}^{\pm}(., t)}(u) \Phi(-u) \frac{du}{2\pi i}
\end{equation}
as an absolutely convergent integral (by Stirling's formula and the decay of $\Phi$), where $0 < \varepsilon < \min_j(- \Re u_j)$   and 
$$ \widehat{\mathcal{J}^{\pm}(., t)}(u) = (2\pi)^{-u} \Gamma(u/2 + it) \Gamma(u/2 - it) \begin{cases} \cos(\pi u/2), & \pm = +,\\ \cosh(\pi t), & \pm = -\end{cases}$$
by  \cite[(3.12)]{BK}.  To deduce \eqref{L1} 
we   shift the contour to the left to $\Re u = -B + 1/2$.  On the way we pick up poles at $u = -2n \pm 2 it$, $n \in \{0, 1, \ldots, [\frac{1}{2}(B-\frac{1}{2})]\}$ with residues of the shape
$$\pm 2 (2\pi)^{2n \pm 2 it} \cosh(\pi t) \Gamma(-n \mp 2it) \Phi(2n \pm 2it) \ll  T^{-2n } (1 + |t|)^{c-B-1/2}$$
with various sign combinations.  
We estimate the remaining integral by Stirling's formula as
\begin{displaymath}
\begin{split}
&\ll T^{-B+1/2} \int_{(-B+1/2)} \big((1 + |\Im u + 2t|)(1 + |\Im u - 2t|)\big)^{-(2B+1)/4} (1 + |u|)^{-B} |du|\\
& \ll T^{c-B + 1/2} (1 + |t|)^{c-B-1/2}.
\end{split}
\end{displaymath}
To obtain analytic continuation to some $t_0$ (distinct from $\pm iu_j/2$) with $|\Im t_0| < \rho$, we shift the contour in \eqref{start} to $\Re u = 2|\Im t_0| + \varepsilon$ for some sufficiently small $\varepsilon > 0$.  We may pick up poles with residue $\widehat{\mathcal{J}^{\pm}(., t)}(-u_j)$, which are meromorphic functions with poles at most at $t = \pm i u_j/2$, while the remaining integral is holomorphic  in a neighbourhood of $t_0$, and \eqref{L2} is clear. 

The holomorphic case is similar, but easier. By \cite[(A.7)]{BLM} we have  $$\widehat{\mathcal{J}^{\text{hol}}(., k)}(u)  = \frac{i^k\pi}{(2\pi)^{u}} \frac{\Gamma( \frac{1}{2}(u + k-1)) }{\Gamma( \frac{1}{2}(1+k-u))},$$  so we can start with the contour $\Re u = 0$ in \eqref{start}, and we have nothing to show for $k \leq B$. For $k \geq B$, the contour shift to the left $\Re u = -  B  +1 $ does not produce any poles. The remaining integral   can be estimated in the same way using the stronger bound $$ \widehat{ \mathcal{J}^{\text{hol}}(., k)}(-B+1 + iw) \ll  \frac{\Gamma((-B + k   + i w)/2)}{\Gamma(  (k + B- iw)/2) } \ll  (k + |w|)^{-B}.$$  \\

Our final lemma will be used to analytically continue the Eisenstein terms $\mathcal{N}^{\text{Eis}}_{q, \ell}(s, w; h)$ and $\mathcal{M}_N^{\text{\rm Eis}}(s, w, \mathfrak{h})$ defined in \eqref{curly-N} and \eqref{curly-M} for $\Re s, \Re w > 1$. 
 
 \begin{lemma}\label{ana-cont} Suppose that $F = {\tt E}_{0}$. 
 
 {\rm (a)} Let $h$ be a $T$-admissible function for some $T \geq 1$. The term $\mathcal{N}^{\text{{\rm Eis}}}_{q, \ell}(s, w; h)$, initially defined in $\Re s, \Re w > 1$ can be continued holomorphically to $\Re s, \Re w > -A$. Its analytic continuation for $\Re s, \Re w < 1$ is given by $\mathcal{N}^{\text{{\rm Eis}}}_{q, \ell}(s, w; h)$ plus the polar term    \begin{equation}\label{pol1}
 \begin{split}
 &\mathcal{R}_{q, \ell}(s, w; h)  :=   \sum \underset{\substack{t =   i(s-1)\\ t =   i(w-1)}}{\text{{\rm res}}}   \frac{ (L(s - it, \chi) \zeta(s + it))^3 L(w - it, \chi) \zeta^{(q)}(w + it)}{L(1 - 2 it, \chi)L(1 + 2 it, \bar{\chi})} \frac{\overline{\Lambda_{\infty, t}}(\ell; w)}{\ell^w} h(t) \\
  \end{split}
 \end{equation}
 {\rm (b)} Let $B \geq 5$ be a constant. Let $\mathfrak{h}_{{\tt s}, {\tt w}} = (h_{{\tt s}, {\tt w}}, h_{{\tt s}, {\tt w}}^{\text{{\rm hol}}})$ where $h_{{\tt s}, {\tt w}}(t)$ is meromorphic in (a neighbourhood of) $|\Im t | \leq 1/2$, $1/2 \leq \Re {\tt s} \leq B$, $1/2 \leq \Re {\tt w} \leq B$  with at most finitely many poles at $t_1({\tt s}, {\tt w}), \ldots, t_n({\tt s}, {\tt w})$ (listed with multiplicity), all of which have $\Im t_j({\tt s}, {\tt w}) \not= 0$. Suppose that $q^2 \nmid N$. The term $\mathcal{M}_N^{\text{\rm Eis}}({\tt s}, {\tt w}, \mathfrak{h}_{{\tt s}, {\tt w}})$,  initially defined in $1 < \Re {\tt s}, \Re {\tt w} < B$ can be continued to (a neighbourhood of) $\Re {\tt s}, \Re {\tt w} \geq 1/2$. Its meromorphic continuation for $1/2 \leq \Re {\tt s}, \Re {\tt w} < 1$ is given by  $\mathcal{M}_N^{\text{\rm Eis}}({\tt s},{\tt  w}, \mathfrak{h}_{{\tt s}, {\tt w}})$ plus the polar term \begin{equation*}
\begin{split} 
 &\widetilde{\mathcal{R}}_{N}({\tt s}, {\tt w}; \mathfrak{h}_{{\tt s}, {\tt w}}) :=     \underset{ t =  \pm  i({\tt s}-1 ) }{\text{{\rm res}}}  (\pm 1) \frac{ \zeta({\tt s}+it)^3\zeta({\tt s}-it)^3  L({\tt w} + it, \bar{\chi} )L({\tt w} - it, \bar{\chi} )}{\zeta(1 + 2 it )\zeta(1 - 2it)}    \tilde{L}_N({\tt s}, {\tt w}; (t, \text{{\rm triv}})) h_{{\tt s}, {\tt w}}(t),
 \end{split}
 \end{equation*}
 which has poles at most at  $\pm i({\tt s}-1 ) = t_{j}({\tt s}, {\tt w})$, $1 \leq j \leq n$, as well as a possible pole at ${\tt s}=1/2$ of order at most $\omega(N/(N, q))$.
 \end{lemma}

\textbf{Proof.} (a) For $t \in \Bbb{R}$ choose a sufficiently small continuous function  $0 < \sigma(t) < 1/4$  such that for $\Re s, \Re w > 2$  the integrand of \eqref{curly-N} is  holomorphic on and in between  the two curves $t \mapsto t \pm i \sigma(t)$.  Let initially $1 <  \Re s < 1 + \sigma(\Im s)$, $1 < \Re w  < 1 + \sigma(\Im w)$ and recall \eqref{eis-expl}. In \eqref{curly-N}  we shift the   $t$-contour up to $\Im t = \sigma(\Re t)$, picking up possible poles at $t = i(s-1), i(w-1)$ in the term corresponding to $\mathfrak{a} = \infty$. The remaining integral and the polar contribution are holomorphic in $1- \sigma(\Im s) <  \Re s < 1 + \sigma(\Im s)$, $1-  \sigma(\Im w) < \Re w  < 1 + \sigma(\Im w)$ and provides continuation to this region. For $s, w$ with $1- \sigma(\Im s) <  \Re s < 1  $, $1-  \sigma(\Im w) < \Re w  < 1  $, we may shift the $t$-contour back to $\Im t = 0$, picking up poles from the $\mathfrak{a} = 0$ term at $t = -i(s-1), -i(w-1)$ which by symmetry just doubles the contribution of the previous poles. This proves the desired formula for the continuation of $\mathcal{N}^{\text{{\rm Eis}}}_{q, \ell}(s, w; h)$ with the correction term as given in \eqref{pol1} for   $1- \sigma(\Im s) <  \Re s < 1  $, $1-  \sigma(\Im w) < \Re w  < 1  $. We now observe that for $\Re s, \Re w < 1$ and $t \in \{i(s-1), i(w-1)\}$, the factor $L(1 + 2it, \bar{\chi})$ in the denominator is zero-free, while the factor $L(1 - 2it, \chi)$ is canceled by the $L$-value in the numerator. Hence both the expression \eqref{pol1} and $\mathcal{N}_{q, \ell}(s, w; h)$ are holomorphic in  $-A < \Re s, \Re w  < 1$ as desired. 

(b) This is the same contour shift argument, and we observe that $q^2 \nmid N$ excludes   characters $\psi$ of conductor $q$, so that $\bar{\chi}\psi$ can never be trivial and hence only the trivial character $\psi = \text{triv}$ in \eqref{curly-M} contributes to the polar term. The meromorphic continuation of  $\tilde{L}_N({\tt s}, {\tt w}; (\pm i({\tt s}-1)), \text{{\rm triv}})) $ to a neighbourhood of $ \Re {\tt s}, \Re {\tt w} > 1/2$ with a possible pole at   ${\tt s} = 1/2$ of order at most $\omega(N/(N, q))$ follows from Lemma \ref{cor-eis} below.

 \section{Functional equations}
 
\subsection{Voronoi on ${\rm GL}(3)$}  We start with a re-formulation of the ${\rm GL}(3)$ Voronoi formula. For
\begin{equation}\label{region1}
\Re  (s + u/2) > 1, \quad \Re(w + u/2) > 1, \quad  \Re u < -1/2
\end{equation}
 we define the absolutely convergent expression
 \begin{equation}\label{defE}
 \mathcal{E}^{\pm}_{q, \ell}(s, u, w) = \sum_{n_1, n_2, c} \sum_{\ell \mid r}    \frac{A(n_2, n_1) \chi(n_1) S (\pm r\bar{q}, n_2, c) \bar{\chi}( r\bar{c})}{n_2^{s+u/2} n_1^{2s} r^{w+u/2  } c^{1-u}}.
\end{equation}
This is well-defined since the character forces $(q, c) = 1$. We also define 
$$\widetilde{\mathcal{E}}^{\pm}_{q, \ell}(s,   u, w) =   \sum_{(c, d) = 1} \sum_{\substack{(m, q) = 1\\ \ell \mid dm}} \sum_n  \frac{A(m, n) \bar{\chi}(d)\chi(c) e(\pm \bar{d} n q/c)}{m^{s+w} n^{1-s-u/2} c^{-1+3s+u/2} d^{w+u/2}}, $$
an expression that is absolutely convergent in 
\begin{equation}\label{abs1}
\Re (s+w) > 1, \quad \Re (w + u/2) > 1, \quad \Re (s + u/2) < 0, \quad \Re(3s + u/2) > 2.
\end{equation}
Let
\begin{equation}\label{defgmu}
\begin{split}
\mathcal{G}_{\mu}^{\pm}(s)   = 4(2\pi)^{-3s} \prod_{j=1}^3 \Gamma(s + \mu_j)\Bigg(\prod_{j=1}^3 \cos\Big(\frac{\pi(s + \mu_j)}{2}\Big) \pm \frac{1}{i}  \prod_{j=1}^3\sin\Big(\frac{\pi(s + \mu_j)}{2}\Big)\Bigg) ,
\end{split}
\end{equation}
where as above $\mu$ denotes the archimedean Langlands parameter of our fixed automorphic automorphic form $F$ for ${\rm SL}_3(\Bbb{Z})$. We recall that $F$ is either cupsidal or the Eisenstein series ${\tt E}_{0}$.

\begin{lemma}\label{func3} Both $\mathcal{E}^{\pm}_{q, \ell}(s,   u, w) $ and $\widetilde{\mathcal{E}}^{\pm}_{q, \ell}(s, u,  w) $ have meromorphic continuation to the region
   \begin{equation}\label{variables}
   \begin{split}
 &\Re (w+u/2)> 1, 
   \quad  \Re(3s + u/2) > 2, \quad \Re (3s - u/2) > 4, \quad  \Re u < -1/2,  \end{split}
 \end{equation}
and satisfy the functional equation 
\begin{equation*} 
\left(\begin{matrix}  \mathcal{E}^{+}_{q, \ell }(s, u, w) \\  \mathcal{E}^{-}_{q, \ell }(s, u, w)\end{matrix}\right)
=  \left(\begin{matrix} \mathcal{G}_{-\mu}^-(1-s-u/2) & \mathcal{G}_{-\mu}^+(1-s-u/2) \\ \mathcal{G}_{-\mu}^+(1-s-u/2) & \mathcal{G}_{-\mu}^-(1-s-u/2) \end{matrix}\right) \left(\begin{matrix}\widetilde{\mathcal{E}}^{+}_{q, \ell }(s, u, w)  \\  \widetilde{\mathcal{E}}^{-}_{q, \ell }(s, u, w)\end{matrix}\right).  \end{equation*}
  If $F$ is cuspidal, both $ \mathcal{E}^{\pm}_{q, \ell }(s, u, w)$ and $\widetilde{\mathcal{E}}^{\pm}_{q, \ell }(s, u, w)$ are holomorphic in \eqref{variables}. If $F = {\tt E}_{0}$, 
the only polar line of $\mathcal{E}^{\pm}_{q, \ell}$ is at $s+u/2 = 1$ (with multiplicity 3), and the  only polar line of $\widetilde{\mathcal{E}}^{\pm}_{q, \ell}$ is at $s+u/2 =  0$ (with multiplicity 3).  We have the upper bound 
 \begin{equation*}
  \begin{split}
 &  \mathcal{E}^{\pm}_{q, \ell}(s, u, w)  \ll_{s, w} (1 + |u|)^{3c( s + \frac{1}{2}u)  +\varepsilon},\\
   & \widetilde{\mathcal{E}}^{\pm}_{q, \ell}(s,  u, w) \ll_{s, w} (1 + |u|)^{3c(1 - s - \frac{1}{2}u)+\varepsilon}
   \end{split}
   \end{equation*}
   with the notation \eqref{cs}, away from the pole if $F = {\tt E}_0$. 
  \end{lemma}
 
 \textbf{Remark:} Note that \eqref{region1} and \eqref{abs1} have no intersection, but are both contained in \eqref{variables}.  \\
 
 \textbf{Proof.}  We  quote from \cite[Section 4]{BK}. For $(c, d) =1$, $\Re v > 1$ let 
\begin{equation*}
\Phi(c, \pm d, m; v) := \sum_{n > 0} A(m, n) e\left(\pm \frac{n \bar{d} }{c}\right) n^{-v},
\end{equation*}
and for $\Re v > 0$ let 
 \begin{equation*}
 \Xi(c, \pm d, m; v) := c \sum_{n_1 \mid cm} \sum_{n_2 > 0} \frac{A(n_2, n_1)}{n_2n_1} S(\pm  md, n_2, mc/n_1) \left(\frac{ n_2n_1^2}{c^3m}\right)^{-v}.
 \end{equation*}
Both functions have meromorphic continuation to $v \in \Bbb{C}$ and satisfy the functional equations \cite[(4.14), (4.15)]{BK}
\begin{equation}\label{func1a}
\left(\begin{matrix}\Phi(c, d, m; v) \\ \Phi(c, -d, m; v)\end{matrix}\right) =  \left(\begin{matrix} \mathcal{G}_{\mu}^+(1-v) & \mathcal{G}_{\mu}^-(1-v) \\ \mathcal{G}_{\mu}^-(1-v) & \mathcal{G}_{\mu}^+(-v) \end{matrix}\right)  \left(\begin{matrix} \Xi(c, d, m; -v) \\ \Xi(c, -d, m; -v)\end{matrix}\right),
\end{equation}
\begin{equation}\label{func1}
\left(\begin{matrix}\Xi(c, d, m; v) \\ \Xi(c, -d, m; v)\end{matrix}\right) =  \left(\begin{matrix} \mathcal{G}_{-\mu}^-(-v) & \mathcal{G}_{-\mu}^+(-v) \\ \mathcal{G}_{-\mu}^+(-v) & \mathcal{G}_{-\mu}^-(-v) \end{matrix}\right)  \left(\begin{matrix} \Phi(c, d, m; -v) \\ \Phi(c, -d, m; -v)\end{matrix}\right).
\end{equation}
Poles can only occur if $F = {\tt E}_0$, in which case $\Phi$ has only a pole at $v = 1$ (of order 3), and $\Xi$ has only a pole (of order 3)  at $v = 0$.

The functions $\Phi$ and $\Psi$ satisfy  the   uniform bounds \cite[(4.16)]{BK}  
\begin{equation}\label{bound}
\begin{split}
& \Phi(c, d, m; v) \ll_{\Re v} \alpha(m)\left(mc^3 (1 + |v|)^3\right)^{ c(v)+\varepsilon},\\
&\Xi(c, d, m; v) \ll_{\Re v} \alpha(m)  (mc^3 )^{ c(-v)+\varepsilon} (1 + |  v|^3)^{c(v+1)+\varepsilon}.
\end{split}
\end{equation}
with $ \alpha(m) =  \max_{d \mid m} |A(d, 1)|$, away from the pole if $F = {\tt E}_0$.  In the region \eqref{variables}  we define 
 \begin{equation*}
 \mathcal{D}^{\pm}_{q, \ell}(s, u, w) := \sum_{(c, d) = 1} \sum_{\substack{(m, q) = 1\\ \ell \mid md}}  \frac{ \bar{\chi}(d)\chi(c)  \Xi(c, \pm d\bar{q}, m; -1 + s + u/2)}{c^{3s+u/2 -1} m^{s+w}d^{w+u/2}}
 \end{equation*}
 and 
 \begin{equation*}
\widetilde{\mathcal{D}}^{\pm}_{q, \ell}(s, u, w) := \sum_{(c, d) = 1} \sum_{\substack{(m, q) = 1\\ \ell \mid md}} \frac{ \bar{\chi}(d)\chi(c) \Phi(c, \pm d\bar{q}, m; 1 - s - u/2)}{c^{3s + u/2 - 1}m^{w+s} d^{w + u/2}}.
\end{equation*} 
(Note that automatically $(c, q) =1$.) It follows from \eqref{bound}   (see \cite[Section 5]{BK})  that both expressions are absolutely convergent in \eqref{variables} and are holomorphic except for  polar lines at $s+u/2 = 1$ of $\mathcal{D}^{\pm}_{q, \ell}$ and  polar lines at $s+u/2 =0 $
of $\widetilde{\mathcal{D}}^{\pm}_{q, \ell}$ if $F = {\tt E}_0$. By \eqref{func1} we have the functional equation
 \begin{equation*} 
\left(\begin{matrix}  \mathcal{D}^{+}_{q, \ell}(s, u, w) \\  \mathcal{D}^{-}_{q, \ell}(s, u, w)\end{matrix}\right)
=  \left(\begin{matrix} \mathcal{G}_{-\mu}^-(-v) & \mathcal{G}_{-\mu}^+(-v) \\ \mathcal{G}_{-\mu}^+(-v) & \mathcal{G}_{-\mu}^-(-v) \end{matrix}\right) \left(\begin{matrix}\widetilde{\mathcal{D}}^{+}_{q, \ell}(s, u, w)  \\ \widetilde{\mathcal{D}}^{-}_{q, \ell}(s, u, w)\end{matrix}\right)  \end{equation*}
as well as the upper bound  \cite[(5.5)]{BK}
 \begin{equation*}
  \begin{split}
 &  \mathcal{D}^{\pm}_{q, \ell}(s, u, w) \ll_{s, w} (1 + |u|)^{3c(s + \frac{1}{2}u) +\varepsilon},\\
   & \widetilde{\mathcal{D}}^{\pm}_{q, \ell}(s, u, w) \ll_{s, w} (1 + |u|)^{3c( 1- s - \frac{1}{2}u) +\varepsilon}
   \end{split}
   \end{equation*}
(away from the possible pole). In the smaller region \eqref{abs1} we have  
$\widetilde{\mathcal{E}}^{\pm}_{q, \ell}(s, u, w) =  \widetilde{\mathcal{D}}^{\pm}_{q, \ell}(s, u, w). $ 
In the smaller region \eqref{region1}  we can recast $\mathcal{D}^{\pm}_{q, \ell}(s, u, w)$ by absolute convergence as  
\begin{displaymath}
\begin{split}
\sum_{(c, d) = 1} \sum_{\substack{(m, q) = 1\\ \ell \mid md}}  \frac{ \bar{\chi}(d)\chi(c) c}{c^{3s+u/2 -1} m^{s+w}d^{w+u/2}}  \sum_{n_1 \mid cm} \sum_{n_2 > 0} \frac{A(n_2, n_1)}{n_2n_1} S(\pm  md\bar{q}, n_2, mc/n_1) \left(\frac{ n_2n_1^2}{c^3m}\right)^{-(s+u/2-1)}. 
\end{split}
\end{displaymath}
The triples $$(m, c, d) \quad \text{with} \quad n_1 \mid cm, (cm, q) = (c, d) = 1, \ell \mid dm$$
are in one-to-one correspondence  with pairs $$(r, C) \quad \text{with} \quad (q, n_1C) = 1, \ell \mid r$$
via
$$C = \frac{cm}{n_1}, \quad r = md,\quad\quad  m = (n_1C, r), \quad d = \frac{r}{m}, \quad c = \frac{n_1C}{m}.$$
Hence
\begin{displaymath}
\begin{split}
\mathcal{D}^{\pm}_{q, \ell}(s, u, w)  = \sum_{(n_1C, q) = 1} \sum_{\ell \mid r}  \sum_{n_2}   \frac{\bar{\chi}(r) \chi(n_1C) A(n_2, n_1) S(\pm r\bar{q}, n_2, C)}{r^{w+u/2}C^{1-u}n_1^{2s}n_2^{s+u/2}}
\end{split}
\end{displaymath}
so that 
$\mathcal{E}^{\pm}_{q, \ell}(s, u, w) =  \mathcal{D}^{\pm}_{q, \ell}(s, u, w). $

\subsection{Laurent coefficients}

We now discuss the Laurent coefficients of  $\mathcal{E}^{\pm}_{q, \ell}$ at $u = 2-2s $ when $F = {\tt E}_0$.   
 Let  
 \begin{equation*}
\mathcal{E}_{q, \ell}^{\pm}(s, u, w) =  \sum_{j=1}^3 \frac{\mathcal{R}_{q,\ell; j}(s, w)}{(u - (2-2s))^j} + O(1). 
\end{equation*}
denote the principal part at $u = 2-2s$. 
We will see in the proof of the next lemma that    $\mathcal{R}_{q,\ell; j}(s, w)$ 
is  independent of the $\pm$ sign. 

\begin{lemma}\label{res} Let $(q, \ell) = 1$. Then  $\mathcal{R}_{q,\ell;j}(s, w)$ is  meromorphic in   the region $$\Re w > 0, \quad \Re s> 1/3.$$  In a neighbourhood of the  region  $1/2 \leq \Re s \leq \Re w$, the  function
$$(s+w-1)^{5-j}\mathcal{R}_{q, \ell; j}(s, w)$$
is holomorphic and bounded by $O_{s, w}( \ell^{-1}(q\ell)^{\varepsilon})$.
\end{lemma}


 \textbf{Proof.}  It is convenient to consider more general Eisenstein series than ${\tt E}_0$. For $\mu \in \Bbb{C}^3$ in a neighbourhood of 0, $\mu_1 + \mu_2 + \mu_3 = 0$ we consider the Eisenstein series ${\tt E}_{\mu}$ with Hecke eigenvalues
 \begin{equation*}
A(n, 1)  = \sum_{abc = n} a^{-\mu_1} b^{-\mu_2} c^{-\mu_3},\quad A(1, n) = \sum_{abc = n} a^{\mu_1} b^{\mu_2} c^{\mu_3}. 
\end{equation*}
 and $A(n_1, n_2)$  given by \eqref{eisen-hecke}. If $\mu_1, \mu_2, \mu_3$ are pairwise distinct, $\mathcal{E}_{q, \ell}^{\pm}(s, u, w)$ has three simple poles at $u = 2-2s-2\mu_j$, $j= 1, 2, 3$, and we compute the residue at $u = 2-2s-2\mu_3$ which we call $\mathcal{R}^{\mu_3}_{q, \ell}(s, w)$. 
  In the region of absolute convergence we have 
$$\mathcal{E}_{q, \ell}^{\pm}(s, u, w)  = 
\sum_{\ell \mid r}  \sum_{ x_1, x_2, x_3,  c} \sum_{d,    b_1, b_2, b_3} \frac{\chi(d) \mu(d)   S(\pm r\bar{q}, db_1b_2b_3, c) \bar{\chi}(r) \chi(x_1x_2x_3c)}{x_1^{2s-\mu_1}x_2^{2s-\mu_2} x_3^{2s-\mu_3} d^{3s + u/2}b_1^{s+u/2 + \mu_1} b_2^{s + u/2+\mu_2} b_3^{s+u/2+\mu_3}c^{1-u} r^{w+u/2}}.$$
We split the $b_3$-sum into residue classes mod $c$, obtaining a Hurwitz zeta-function $$\sum_{b_3 \equiv \beta\, (\text{mod } c)} \frac{1}{b_3^{s+u/2 + \mu_3}} = \frac{1/c}{(s + u/2 + \mu_3 - 1)} + O(1).$$
For reasons of convergences let us initially assume $\Re \mu_2, \Re \mu_1 > \Re \mu_3$. 
Then we see that   $\mathcal{R}^{\mu_3}_{q, \ell}(s, w)$ equals
\begin{displaymath}
\begin{split}
  & \sum_{\ell \mid r}  \sum_{ x_1, x_2, x_3,  c} \sum_{d,    b_1, b_2 }  \sum_{b_3 \, (\text{mod } c)} \frac{ \chi(d)\mu(d)   S(\pm r\bar{q}, db_1b_2b_3, c) \bar{\chi}(r) \chi(x_1x_2x_3c)}{x_1^{2s-\mu_1}x_2^{2s-\mu_2} x_3^{2s-\mu_3} d^{1+2s-\mu_3}b_1^{1 + \mu_1 - \mu_3}b_2^{1 + \mu_2 - \mu_3}    c^{2s  + 2\mu_3} r^{1-s+w- \mu_3}}\\
& = \prod_{j=1}^3 L(2s-\mu_j, \chi) \sum_{\ell \mid r}  \sum_{    c} \sum_{c \mid d    b_1 b_2 }   \frac{ \chi(d)\mu(d)   S(\pm r\bar{q}, 0, c) \bar{\chi}(r) \chi( c)}{  d^{1+2s-\mu_3}b_1^{1 + \mu_1 - \mu_3}b_2^{1 + \mu_2 - \mu_3}    c^{2s -1 + 2\mu_3} r^{1-s+w- \mu_3}}\\
& = \prod_{j=1}^3 L(2s-\mu_j, \chi) \sum_{  c} \sum_{a \mid c} \sum_{\ell \mid ra}  \sum_{c \mid d    b_1 b_2 }   \frac{ \chi(d) a\mu(c/a)\mu(d)    \bar{\chi}(ra) \chi( c)}{ d^{1+2s-\mu_3}b_1^{1 + \mu_1 - \mu_3}b_2^{1 + \mu_2 - \mu_3}    c^{2s -1 + 2\mu_3} (ra)^{1-s+w- \mu_3}}.
\end{split}
\end{displaymath}
In particular we see the sign-independence. 
The quadruple sum is an Euler product whose $p$-Euler factors can be computed by brute force using geometric series, for instance by distinguishing the cases
\begin{itemize}
\item $p \nmid d$,  $v_p(b_1) \geq v_p(c) + 1$;
\item  $p \nmid d$,  $v_p(b_1) < v_p(c)$, $v_p(b_2) \geq v_p(c) - v_p(b_1)$;
\item $v_p(d) = 1$, $v_p(c) = v_p(a) + 1$, $v_p(b_1) \geq v_p(a) + 1$;
\item $v_p(d) = 1$, $v_p(c) = v_p(a) + 1$, $v_p(b_1) \leq v_p(a)$, $v_p(b_2) \geq v_p(a) - v_p(b_1)$;
\item  $v_p(d) = 1$, $v_p(c) = v_p(a) \geq 1$, $v_p(b_1) \geq v_p(a)$;
\item  $v_p(d) = 1$, $v_p(c) = v_p(a) \geq 1$, $v_p(b_1) \leq v_p(a)-1$, $v_p(b_2) \geq v_p(a) - v_p(b_1) - 1$;
\item $v_p(d) = 1$, $p \nmid ac$. 
\end{itemize}
 In this way we see that for $p \nmid \ell$ the $p$-Euler factor equals (using that $\mu_1 + \mu_2 + \mu_3 = 0$)
\begin{displaymath}
\begin{split}
\frac{P_p(s, w; \mu_1, \mu_2, \mu_3)}{(1 - p^{-1+\mu_3 - \mu_1}) (1 - p^{-1+\mu_3 - \mu_2})  (1 - \chi_0(p) p^{-s-w- \mu_1}) (1 - \chi_0(p) p^{-s-w- \mu_2}) (1 - \bar{\chi}(p) p^{-1-w+s+\mu_3})}
\end{split}
\end{displaymath}
where $P_p(s, w; \mu_1, \mu_2, \mu_3)$ is defined as
$$1   - \frac{\chi(p)}{p^{2s -\mu_1}} - \frac{\chi(p)}{p^{2s-\mu_2}}  + \frac{\chi(p)}{p^{1+3s+w-\mu_3 + \mu_1}}+ \frac{\chi(p)}{p^{1+3s+w-\mu_3 + \mu_2}} - \frac{\chi_0(p)}{p^{1+s+w - \mu_3 + \mu_1+\mu_2}}   + \frac{\chi(p)^2}{p^{4s+\mu_3}}- \frac{\chi(p)^2}{p^{1+5s+w -\mu_3}}$$
where $\chi_0$ is the trivial character modulo $q$. Note that
$$
\prod_p P_p(s, w; \mu_1, \mu_2, \mu_3) = \frac{P(s, w; \mu_1, \mu_2, \mu_3)}{L(s - 2\mu_1, \chi)L(s - 2\mu_2, \chi)}$$
where $P(s, w; \mu_1, \mu_2, \mu_3)$ is holomorphic  and absolutely bounded in a neighbourhood  of $\Re s \geq 1/2$, $\Re w \geq 0$, $\mu = 0$. 
The $p$-Euler factor for $p \mid \ell$   converges absolutely and can be bounded trivially by
$$O\left(p^{-v_p(\ell) \Re (s + w)} + p^{-v_p(\ell)\Re(1+w-s)}\right).$$
 Hence $\mathcal{R}^{\mu_3}_{q, \ell}(s, w)$ equals
\begin{displaymath}
\begin{split}
 &   \zeta(1+\mu_1 - \mu_3)\zeta(1+\mu_2 - \mu_3) \zeta^{(q)}(s+w+\mu_1)\zeta^{(q)}(s+w+\mu_2) L(2s-\mu_3, \chi)  L(1+w-s-\mu_3, \bar{\chi}) \mathcal{L}_{\mu}(s, w)
  \end{split}
  \end{displaymath}
where $\mathcal{L}_{\mu}(s, w)$ is holomorphic  in  a neighbourhood of $\Re s \geq 1/2$, $\Re w \geq 0$, $\mu = (\mu_1, \mu_2, \mu_3) = 0$ and bounded by
$$O\left((\ell^{- \Re (s + w)   }+ \ell^{- \Re(1+w-s) })(\ell q)^{\varepsilon}\right).$$
Similar formulae with exchanged indices hold for 
$\mathcal{R}^{\mu_1}_{q, \ell;1}(s, w) $ and $\mathcal{R}^{\mu_2}_{q, \ell;1}(s, w) $. Taking   limits $\mu_1, \mu_2, \mu_3 \rightarrow 0$ (on the hyperplane $\mu_1 + \mu_2 + \mu_3 = 0$), we can obtain expressions for $\mathcal{R}_{q, \ell; j}(s, w)$, $j = 1, 2, 3$, that involve derivatives up to order $3-j$ of 
$$ \zeta^{(q)}(s+w)^2  L(2s, \chi)  L(1+w-s, \bar{\chi}) \mathcal{L}_{(0, 0, 0)}(s, w).$$
The lemma follows.

\subsection{Poisson summation}

In this subsection we present a re-formulation of Poisson summation in terms of the functional equation of the Hurwitz zeta-function. For $\Re v > 0$ we define the two absolutely convergent expressions
\begin{equation}\label{defTheta}
\Theta_q(n, d; v) = \sum_{(c, d) = 1} \frac{\chi(c) e(  \bar{c} n q/d)}{c^{v}}, \quad \tilde{\Theta}(n, d; v) =  \sum_c \frac{  \bar{\chi}(c) \chi(  d) S(n,   c, d)}{c^{v}}.
\end{equation}
Let
 \begin{equation}\label{defG}
 G^{\pm}(s) =  \Gamma(s) (2\pi)^{-s} \exp(\pm i\pi s/2). 
 \end{equation}

\begin{lemma}\label{func-gl1} For $(d,q) = 1$,  the functions  $\Theta_q(n, d; v) $  and $\tilde{\Theta}(n, d; v)$ have holomorphic  continuation to all of $\Bbb{C}$ and satisfy the functional equation
$$\left(\begin{matrix} \Theta_q(n, d; v)\\ \Theta_q(-n, d; v) \end{matrix}\right)    = \frac{\tau(\chi)}{(qd)^v}  \left(\begin{matrix} G^-(1-v) & G^+(1-v)\\ G^+(1-v) & G^-(1-v) \end{matrix}\right) \left(\begin{matrix} \tilde{\Theta}(n, d; 1-v)\\ \tilde{\Theta}(-n, d; 1-v) \end{matrix}\right)  $$
as well as the uniform upper bounds
\begin{equation}\label{theta}
\begin{split}
 &\Theta_q(n, d; v) \ll_q (d(1+|v|))^{c(v)+\varepsilon}\\
 &\tilde{\Theta}(n, d; v) \ll_q d^{1/2} (d(1+|v|))^{c(v)+\varepsilon}
 \end{split}
 \end{equation}
 \end{lemma}

\textbf{Proof.} Using the Hurwitz zeta function, we have 
$$\Theta_q(n, d; v) = \frac{1}{(qd)^v}\underset{ \gamma \, (\text{mod }dq)}{\left.{\sum}\right.^{\ast}} \chi(\gamma)e\left( \frac{\bar{\gamma} nq}{ d}\right) \zeta\left(v, \frac{\gamma}{qd}\right)$$
so that meromorphic continuation of $\Theta_q(n, d; v)$ follows from properties of the Hurwitz zeta-function. For $\Re v \geq  1+\varepsilon$ we have the uniform bound 
$\Theta_q(n, d, v) \ll 1$. The residue at the only possible pole $v=1$ equals
$$\frac{1}{qd}
\underset{ \gamma \, (\text{mod }dq)}{\left.{\sum}\right.^{\ast}}  \chi(\gamma)e\left( \frac{\bar{\gamma} nq}{ d}\right) = 0$$ 
(since $(d, q) = 1$ and $\chi$ is non-trivial). The functional equational for the Hurwitz zeta-function \cite[(4.7)]{BH} implies
$$\Theta_q(n, d; v) =\frac{1}{(qd)^v}\underset{ \gamma \, (\text{mod }dq)}{\left.{\sum}\right.^{\ast}} \chi(\gamma)e\left( \frac{\bar{\gamma} nq}{ d}\right) \sum_{\pm} G^{\mp}(1-v) \zeta^{\ast}\left(v, \pm\frac{\gamma}{qd}\right) $$
where
$$\zeta^{\ast}(v, \alpha) = \sum_{c> 0} e(c \alpha) c^{-v}$$
for $\Re v > 1$. For $(d, q) = 1$ we have
\begin{displaymath}
\begin{split}
\underset{ \gamma \, (\text{mod }dq)}{\left.{\sum}\right.^{\ast}}  \chi(\gamma)e\left( \frac{\bar{\gamma} nq}{ d}\right)e\left(\pm \frac{c\gamma}{qd}\right) = \tau(\chi) \bar{\chi}(c) \chi(  d) S(\pm n,   c, d). 
\end{split}
\end{displaymath}
Hence in $\Re v < 0$ we have
$$\Theta_q(n, d, v) = \frac{\tau(\chi)}{(qd)^v} \sum_{\pm} G^{\mp}(1-v) \sum_c \frac{  \bar{\chi}(c) \chi(  d)S(\pm n, c, d)}{c^{1-v}} \ll (d(1+|v|))^{1/2 -\Re v},$$
and the functional equation is clear. The upper bound for $\Theta_q(n, d; v)$ follows from the Phragm\'en-Lindel\"of principle. 
In $\Re v > 1$ we have
$$\tilde{\Theta}(n, d, v) = \frac{1}{(dq)^v} \sum_{\gamma \, (\text{mod } dq)}  \bar{\chi}(\gamma) \chi(  d) S(n, \gamma, d) \zeta\left(v, \frac{\gamma}{dq}\right)$$
and the only possible pole at $v=1$ has residue
$$ \sum_{\gamma \, (\text{mod }dq)} \chi(\bar{\gamma}d) S(n, \gamma, d) = 0$$
(since $(d, q) = 1$ and $\chi$ is non-trivial). Hence also $\tilde{\Theta}(n, d; v)$ can be continued holomorphically to all of $\Bbb{C}$. Inverting the functional equation by 
$$\left(\begin{matrix} G^-(1-v) & G^+(1-v)\\ G^+(1-v) & G^-(1-v) \end{matrix}\right) ^{-1} = \left(\begin{matrix} G^+(v) & G^-(v)\\ G^-(v) & G^+(v) \end{matrix}\right),$$
we combine the trivial bound $\tilde{\Theta}(n, d; v)  \ll d^{1/2+\varepsilon}$ in $\Re v \geq 1+\varepsilon$ with the bound $\tilde{\Theta}(n, d, v) \ll d^{1-v} (1+|v|)^{1/2-\Re v}$ in $\Re v \leq -\varepsilon$. \\
 
 We now define 
\begin{equation}\label{defAB}
\begin{split}
& \mathcal{A}^{\pm}(s, u, w) := \sum_{\substack{(m, q) = 1\\ \ell \mid dm}} \sum_n  \frac{A(m, n) \bar{\chi}(d) \Theta_q(\pm  n, d, -1+3s+u/2 )}{m^{s+w} n^{1-s-u/2}   d^{w+u/2}},\\
&
 \mathcal{B}^{\pm}(s, u, w) := \sum_{\substack{(m, q) = 1\\ \ell \mid dm}} \sum_n  \frac{A(m, n) \bar{\chi}(d) \tilde{\Theta}( \pm  n, d, 2-3s-u/2 )}{m^{s+w} n^{1-s-u/2}   d^{3s+u+w-1}}.
 \end{split}
 \end{equation}
By \eqref{theta}, both of these are absolutely convergent in 
\begin{equation}\label{newvariables}
\begin{split}
&\Re(s+w) > 1, \quad \Re(1 - s- u/2) > 1, \\
& \Re(w+u/2-\xi) > 1, \quad \Re (3s+u+w) > 5/2, \quad \Re(3s/2 + 3u/4 + w) > 2,
\end{split}
\end{equation}
satisfy the functional equation
\begin{equation}\label{func-new}
\left(\begin{matrix} \mathcal{A}^{+}(s, u, w)\\ \mathcal{A}^{-}(s, u, w) \end{matrix}\right)    =   \frac{\tau(\chi)}{q^{-1+3s+u/2 }} \left(\begin{matrix} G^-(2-3s-u/2) & G^+(2-3s-u/2 )\\ G^+(2-3s-u/2 ) & G^-(2-3s-u/2 ) \end{matrix}\right) \left(\begin{matrix} \mathcal{B}^{+}(s, u, w)\\ \mathcal{B}^{-}(s, u, w) \end{matrix}\right),  
\end{equation}
and are bounded by 
\begin{equation}\label{414a}
\begin{split}
& \mathcal{A}^{\pm}(s, u, w)  \ll (1+|u |)^{c(-1+3s+u/2 )+ \varepsilon},\\
&\mathcal{B}^{\pm}(s, u, w)  \ll (1+|u|)^{c(2-3s-u/2) + \varepsilon}.
\end{split}
\end{equation}
We can open $\mathcal{A}^{\pm}(s, u, w)$ in the subregion \eqref{abs1} 
of \eqref{newvariables}, and we can open $\mathcal{B}^{\pm}(s, u, w)$ in the subregion
\begin{equation}\label{b}
 \Re(s+w) > 1, \quad \Re(1 - s- u/2) > 1, \quad \Re (3s+u+w) > 5/2, \quad \Re(2 - 3s - u/2) > 1.
\end{equation}

 \section{Local factors}
 
\subsection{The local Rankin-Selberg factor} For a prime $p$ let $\alpha_{f, \nu}(p)$ ($\nu = 1, 2$), $\alpha_{F, j}(p)$ ($j = 1, 2, 3$) denote the Satake parameters of $f \in \mathcal{B}(q, \chi)$ and $F$ at   $p$ satisfying
\begin{equation*}
 \alpha_{F, 1}(p) \alpha_{F, 2}(p) \alpha_{F, 3}(p)
 = 1, \quad \alpha_{f, 1}(p) \alpha_{f, 2}(p) = \chi(p). 
 \end{equation*}
 We have
\begin{equation*}
\lambda_f(p^{\nu}) = \frac{\alpha_{f, 1}(p)^{\nu+1}  - \alpha_{f, 2}(p)^{\nu+1} }{\alpha_{f, 1}(p) - \alpha_{f, 2}(p)}
\end{equation*}
and 
\begin{equation*}
 A(p^{\nu}, p^{\mu}) =  \det\left(\begin{matrix} \alpha_{F, 1}(p)^{\nu+\mu+2} & \alpha_{F, 2}(p)^{\nu+\mu+2} & \alpha_{F, 3}(p)^{\nu+\mu+2}\\ \alpha_{F, 1}(p)^{\mu+1} & \alpha_{F, 2}(p)^{\mu+1} & \alpha_{F, 3}(p)^{\mu+1}\\1 & 1 & 1\end{matrix}\right) V_F(p)^{-1}
\end{equation*}
where 
$V_F(p)= \det\left(\begin{matrix} \alpha_{F, 1}(p)^{ 2} & \alpha_{F, 2}(p)^{ 2} & \alpha_{F, 3}(p)^{2}\\ \alpha_{F, 1}(p) & \alpha_{F, 2}(p) & \alpha_{F, 3}(p)\\1 & 1 & 1\end{matrix}\right). 
$
Let $$L_p(s, f \times F) 
= \prod_{j=1}^3\prod_{\nu=1}^2 \left(1 - \frac{\alpha_{F, j}(p) \alpha_{f, \nu}(p)}{p^s}\right)^{-1}$$
denote the local Rankin-Selberg factor at $p$. A straightforward computation with geometric series as in \cite[Lemma 13]{BK} shows
\begin{equation}\label{local-RS}
L_p(s, f \times F) = \sum_{n_1, n_2 \mid p^{\infty}} \frac{A(n_2, n_1) \chi(n_1) \lambda_f(n_2)}{n_2^s n_1^{2s}} 
\end{equation}
 for all primes $p$. This formula holds for any character $\chi$ modulo $q$, including the trivial character. The Euler factor is absolutely bounded from above and below in $\Re s > \vartheta + \theta$ ($\leq 5/14 + 7/64 < 1/2$). The same formula holds for $\mathcal{E}_{\mathfrak{a}}(t)$ and $\lambda_{\mathfrak{a}, t}(n_2)$ in place of $f$ and $\lambda_f(n_2)$. 
 
\subsection{The cuspidal correction factors} We start by considering
$$
\sum_{\ell_1\ell_2 = \ell} \frac{1}{\ell_1^{2s}}  \sum_{(m, \ell_2 q) = 1}\sum_{n, c}  \frac{ \bar{\chi}(c)  A(\ell_1m, n)   }{c^{w } m^{2s} n^{s }   } \mathcal{A}^{\text{Maa{\ss}}}_{\ell_2Q}(\pm  n, c; \mathfrak{h}) $$
for $Q \in \{1, q\}$,  $\mathfrak{h} = (h, h^{\text{hol}})$ and  $\Re s , \Re w > 1$ with the notation as in \eqref{defA}. We write this as
$$\sum_{\ell_0 \mid \ell Q} \sum_{f \in \mathcal{B}^{\ast}(\ell_0, \text{triv})}\epsilon_f^{(1 \mp 1)/2} h(t_f) \mathcal{S}(s, w; f)$$
where
$$\mathcal{S}(s, w; f) = \sum_{\ell_1\ell_2 = \ell} \frac{1}{\ell_1^{2s}} \sum_{M \mid \frac{\ell_2Q}{\ell_0}}  \sum_{(m, \ell_2 q) = 1}\sum_{n, c}  \frac{\rho_{f, M, \ell_2Q}(n) \rho_{f, M, \ell_2Q}(c)  \bar{\chi}(c)A(\ell_1m, n)}{c^{w } m^{2s} n^{s'} } $$
By \eqref{rho-cusp},  we have
\begin{displaymath}
\begin{split}
&L^{\ast}(1, \text{Ad}^2 f) \prod_{p \mid \ell_0} (1 - p^{-2})^{-1} \mathcal{S}(s, w; f) \\
&=  \sum_{\ell_1\ell_2 = \ell} \frac{1}{\ell_1^{2s}} \sum_{\ell_0M \mid \ell_2Q}  \sum_{(m, \ell_2 q) = 1} \sum_{\delta_1, \delta_2 \mid M} \sum_{n, c} \frac{\xi_f(M, \delta_1) \xi_f(M, \delta_2)}{Q\ell_2 \nu(Q\ell_2)} \frac{\delta_1\delta_2}{M}  \frac{\lambda_{f}(n/\delta_1) \lambda_{f}(c/\delta_2) \bar{\chi}(c)A(\ell_1m, n)}{c^{w } m^{2s} n^{s } } \\
&  = \sum_{\ell_1\ell_2 = \ell} \frac{1}{\ell_1^{2s}} \sum_{\ell_0M \mid \ell_2Q}  \sum_{(m, \ell_2 q) = 1} \sum_{\delta_1, \delta_2 \mid M} \sum_{n} \frac{\xi_f(M, \delta_1) \xi_f(M, \delta_2)}{Q\ell_2 \nu(Q\ell_2)} \frac{\delta_1^{1-s}\delta_2^{1-w} \bar{\chi}(\delta_2)}{M} \\
& \quad\quad\quad\quad \frac{\lambda_{f}(n ) A(\ell_1m, \delta_1n)}{  m^{2s} n^{s } } L(w, f \times \bar{\chi}). 
\end{split}
\end{displaymath}
We   recognize the $n, m$-sum as $L(s, f \times \tilde{F})$ up to Euler factors at primes dividing $q \ell$ and obtain
$$\mathcal{S}(s, w; f) = \frac{L(s, f \times \tilde{F})  L(w, f \times \bar{\chi})}{L(1, \text{Ad}^2 f)}  \tilde{L}_{\ell Q}(s, w; f) $$
where
\begin{equation}\label{cusp-corr}
\begin{split}
 \tilde{L}_{\ell Q}(s, w; f) = &\frac{L(1, \text{Ad}^2 f)}{L^{\ast}(1, \text{Ad}^2 f)} \prod_{p \mid \ell_0} (1 - p^{-2})  \prod_{p \mid q\ell}L_p(s, f \times \tilde{F}) ^{-1} \sum_{\ell_1\ell_2 = \ell} \frac{1}{\ell_1^{2s}} \sum_{\ell_0M \mid \ell_2Q} \sum_{\delta_1, \delta_2 \mid M}  \\
 &\underset{n, m \mid (\ell q)^{\infty}}{\sum_{(m, \ell_2 q) = 1}  \sum_{n}} \frac{\xi_f(M, \delta_1) \xi_f(M, \delta_2)}{Q\ell_2 \nu(Q\ell_2)} \frac{\delta_1^{1-s}\delta_2^{1-w}\bar{\chi}(\delta_2)}{M}  \frac{\lambda_{f}(n ) A(\ell_1m, \delta_1n)}{  m^{2s} n^{s } } .
 \end{split}
 \end{equation}
Although the right hand side depends on $\ell$ and $Q$ individually, this defines $ \tilde{L}_{N}(s, w; f)$ for every positive integer with $q^2 \nmid N$ by decomposing uniquely $N = \ell Q$ with $(\ell, q ) = 1$, $Q \in \{1, q\}$. 
 In $\Re s, \Re w \geq 1/2$ we use \eqref{xi-arithmetic} and estimate trivially
 \begin{equation}\label{corr-bound}
 \tilde{L}_{\ell Q}(s, w; f) \ll (q\ell)^{\varepsilon}\sum_{\ell_1\ell_2 = \ell} \frac{1}{\ell_1} \sum_{\ell_0M \mid \ell_2Q} \sum_{\delta_1, \delta_2 \mid M} \left(\frac{M^2}{\delta_1\delta_2}\right)^{\vartheta} \frac{(\delta_1\delta_2)^{1/2} (\ell_1\delta_1)^{\theta}}{Q\ell_2 M} \ll (Q\ell)^{\theta - 1} (q\ell)^{\varepsilon}
 \end{equation}
 confirming \eqref{corr-bound-final} in the cuspidal case. With this notation we have
 \begin{equation}\label{AtoM1}
 \sum_{\ell_1\ell_2 = \ell} \frac{1}{\ell_1^{2s}}  \sum_{(m, \ell_2 q) = 1}\sum_{n, c}  \frac{ \bar{\chi}(c)  A(\ell_1m, n)   }{c^{w } m^{2s} n^{s }   } \mathcal{A}^{\text{Maa{\ss}}}_{\ell_2Q}(\pm  n, c; \mathfrak{h})  = \mathcal{M}^{\text{Maa{\ss}}, \pm}_{\ell Q} (s, w; \mathfrak{h}).
 \end{equation}
The same analysis holds for holomorphic cusp  forms. 

\subsection{The Eisenstein correction factors} Similarly as in the preceding subsection we consider now 
$$
\sum_{\ell_1\ell_2 = \ell} \frac{1}{\ell_1^{2s}}  \sum_{(m, \ell_2 q) = 1}\sum_{n, c}  \frac{ \bar{\chi}(c)  A(\ell_1m, n)   }{c^{w } m^{2s} n^{s }   } \mathcal{A}^{\text{Eis}}_{\ell_2Q}(\pm  n, c; \mathfrak{h})  = \sum_{d_0^2 \mid \ell Q} \sum_{\substack{\psi \, (\text{mod } d_0)\\ \text{primitive}}} \int_{\Bbb{R}} \mathcal{S}(s, w; (t, \psi)) h(t) \frac{dt}{2\pi}$$
for $Q \in \{1, q\}$ 
 where $\mathcal{S}(s, w; (t, \psi))$ is defined as
\begin{displaymath}
\begin{split}
&\sum_{\ell_1\ell_2 = \ell} \frac{1}{\ell_1^{2s}}  \sum_{d_0 \mid M_1 \mid d_0^{\infty}} \sum_{\substack{ (M_2, d_0) = 1 \\ d_0 M_1M_2 \mid \ell_2Q}} \sum_{(m, \ell_2 q) = 1}\sum_{n, c}     \frac{  \bar{\chi}(c)  A(\ell_1m, n)  \rho_{\psi, d_0M_1M_2, \ell_2Q}(n, t) \overline{\rho_{\psi, d_0M_1M_2, \ell_2Q}(c, t)}}{c^{w } m^{2s} n^{s } }
\end{split}
\end{displaymath}
for $t\in \Bbb{R}$ and $\psi$ a primitive Dirichlet character modulo $d_0$. We define 
\begin{equation}\label{lchitd2}
L(\psi, t, \ell_2) := \prod_{p \mid \ell_2} \left(\left(1 - \frac{\psi^2(p)}{p^{1 + 2it}}\right) \left(1 - \frac{\bar{\psi}^2(p)}{p^{1 - 2it}}\right)\right)^{-1},
\end{equation}
insert \eqref{rho-eis} and recast $  |L(1 + 2it, \psi^2)|^2\mathcal{S}(s, w; (t, \psi))$ as
\begin{displaymath}
\begin{split}
& \sum_{\ell_1\ell_2 = \ell} \frac{L(\psi, t, \ell_2) }{\ell_1^{2s} \ell_2Q\nu(\ell_2Q)}  \sum_{d_0 \mid M_1 \mid d_0^{\infty}} \sum_{\substack{ (M_2, d_0) = 1 \\ d_0 M_1M_2 \mid \ell_2Q}} \sum_{\delta_1, \delta_2 \mid M_2} \frac{M_1\delta_1\delta_2 \mu(M_2/\delta_1) \mu(M_2/\delta_2)\bar{\psi}(\delta_1) \psi(\delta_2)   }
 { \tilde{\mathfrak{n}}(d_0M_1M_2)^2M_2}  \\
 &  \sum_{(m, \ell_2q) = 1}   \sum_{\substack{c_1,f_1\\ (c_1, \frac{\ell_2Q}{d_0M_1M_2}) = 1}}   \sum_{\substack{  c_2,f_2\\ (c_2, \frac{\ell_2Q}{d_0M_1M_2}) = 1}}    \frac{\bar{\chi}(c_2f_2)A(\ell_1m, c_1f_1M_1\delta_1)\bar{\psi}(c_1f_2)  \psi(c_2f_1) (\frac{c_2f_2M_1\delta_2}{  c_1f_1M_1\delta_1})^{it}  }{(c_2/c_1)^{2it}(c_1f_1M_1\delta_1)^{s  }m^{2s} (c_2f_2 M_1 \delta_2)^{w}}.
\end{split}
\end{displaymath}
The second line can be simplified as
$$ L^{(N)}(w + it, \psi \bar{\chi}) L^{(N)}(w -it, \bar{\psi}\bar{\chi}) \sum_{(m, \ell_2q) = 1}   \sum_{\substack{c_1,f_1\\ (c_1, N) = 1}}      \frac{ A(\ell_1m, c_1f_1M_1\delta_1)\bar{\psi}(c_1 )  \psi( f_1)    }{ c_1^{-2it}(c_1f_1M_1\delta_1)^{s +it }m^{2s} (  M_1 \delta_2)^{w-it}}$$
with $N = \ell_2Q/(d_0M_1M_2)$.  
In the $m, c_1, f_1$-sum we recognize $L(s + it, \tilde{F} \times \psi) L(s - it, \tilde{F} \times \bar{\psi})$ up to Euler factors at primes dividing $\ell q$. Thus, by brute force, we write
$$\mathcal{S}(s, w; (t, \psi)) = \frac{L(s + it, \tilde{F} \times \psi) L(s - it, \tilde{F} \times \bar{\psi})L (w + it, \psi \bar{\chi}) L (w -it, \bar{\psi}\bar{\chi}) }{ |L(1 + 2it, \psi^2)|^2} \tilde{L}_{\ell Q}(s, w; (t, \psi))$$
where $\tilde{L}_{\ell Q}(s, w; (t, \psi))$ for a primitive character $\psi$ of conductor $d_0$ is defined as
\begin{equation}\label{Eis-corr}
\begin{split}
 & \sum_{\ell_1\ell_2 = \ell} \frac{L(\psi, t, \ell_2) }{\ell_1^{2s} \ell_2Q\nu(\ell_2Q)}  \sum_{d_0 \mid M_1 \mid d_0^{\infty}} \sum_{\substack{ (M_2, d_0) = 1 \\ d_0 M_1M_2 \mid \ell_2Q}} \sum_{\delta_1, \delta_2 \mid M_2} \frac{M_1\delta_1\delta_2 \mu(M_2/\delta_1) \mu(M_2/\delta_2)\bar{\psi}(\delta_1) \psi(\delta_2)   }
 { \tilde{\mathfrak{n}}(d_0M_1M_2)^2M_2} \\
 & \sum_{\substack{(m, \ell_2q) = 1\\ m \mid \ell^{\infty}}}   \sum_{\substack{c_1f_1\mid (\ell q)^{\infty} \\ (c_1, \frac{\ell_2Q}{d_0M_1M_2}) = 1}}      \frac{ A(\ell_1m, c_1f_1M_1\delta_1)\bar{\psi}(c_1 )  \psi( f_1)    }{ c_1^{-2it}(c_1f_1M_1\delta_1)^{s +it }m^{2s} (  M_1 \delta_2)^{w-it}}\\
& \prod_{p \mid \frac{\ell_2Q}{d_0M_1M_2}} \left(1 - \frac{\psi \bar{\chi}(p)}{p^{w+it}}\right)\left(1 - \frac{\overline{\psi \chi}(p)}{p^{w-it}}\right) \prod_{p \mid \ell q} \Big(L_p(s + it, \tilde{F} \times \psi)  L_p(s - it, \tilde{F} \times \bar{\psi})\Big)^{-1}
\end{split}
\end{equation} 
with the notation \eqref{lchitd2}.  As in \eqref{cusp-corr}, the formula \eqref{Eis-corr} defines $\tilde{L}_{N}(s, w; (t, \psi))$ for every $N \in \Bbb{N}$ with $q^2 \nmid N$. 
With this notation we have 
 \begin{equation}\label{AtoM2}
 \sum_{\ell_1\ell_2 = \ell} \frac{1}{\ell_1^{2s}}  \sum_{(m, \ell_2 q) = 1}\sum_{n, c}  \frac{  \bar{\chi}(c)  A(\ell_1m, n)   }{c^{w } m^{2s} n^{s }   } \mathcal{A}^{\text{Eis}}_{\ell_2Q}(\pm  n, c; \mathfrak{h}) = \mathcal{M}^{\text{Eis}}_{\ell Q} (s, w; \mathfrak{h}).
 \end{equation}
For $\Re s, \Re w \geq 1/2$, $t \in \Bbb{R}$ we estimate trivially
$$\tilde{L}_{\ell Q}(s, w; (t, \psi)) \ll (\ell q)^{\varepsilon}  \sum_{\ell_1\ell_2 = \ell} \frac{1}{\ell Q}  \sum_{d_0 \mid M_1 \mid d_0^{\infty}} \sum_{ d_0 M_1M_2 \mid \ell_2Q} \sum_{\delta_1, \delta_2 \mid M_2} \frac{(\delta_1\delta_2)^{1/2}}{M_2}(\ell_1M_1\delta_1)^{\theta} \ll (\ell Q)^{\theta-1} (\ell q)^{\varepsilon},$$
 confirming \eqref{corr-bound-final} in the Eisenstein case. In the special case $\psi = \text{triv}$ of conductor $d_0 = 1$ and $F = {\tt E}_{0}$, the following lemma provides analytic continuation of $\tilde{L}_{\ell Q}(s, w; (t, \psi))$ to certain complex values of $t$. 
 
 \begin{lemma}\label{cor-eis} Let $F = {\tt E}_{0}$ and $N \in \Bbb{N}$ with $q^2 \nmid N$. The functions $\tilde{L}_N(s, w; (\pm i(s-1) , \text{{\rm triv}})) $, initially defined in $\Re s, \Re w > 1$ as absolutely convergent series, have meromorphic continuation to an $\varepsilon$-neighbourhood of $\Re s, \Re w \geq 1/2$ with polar divisors at most at $s = 1/2$  and satisfy the bounds
\begin{equation}\label{bound-cor-eis}
 (s - 1/2)^{\omega(\ell)}\tilde{L}_N(s, w; (\pm i(s-1 ), \text{{\rm triv}}))  \ll_{s, w} N^{-1} (qN)^{ \varepsilon }
 \end{equation}
(where $\ell = N/(N, q)$) for $1/2 -\varepsilon \leq \Re s, \Re w < 1$. 
\end{lemma}

\textbf{Proof.} A small variation of \cite[Lemma 14]{BK} (that includes the characters $\bar{\psi}(c_1 )  \psi( f_1)  $) shows that for $d_0M_1M_2 \mid \ell_2 Q$, $Q \in \{1, q\}$, $\ell_1\ell_2 = \ell$ the triple sum
$$\sum_{\substack{(m, \ell_2q) = 1\\ m \mid \ell^{\infty}}}   \sum_{\substack{c_1f_1\mid (\ell q)^{\infty} \\ (c_1, \frac{\ell_2Q}{d_0M_1M_2}) = 1}}      \frac{ A(\ell_1m, c_1f_1M_1\delta_1)\bar{\psi}(c_1 )  \psi( f_1)    }{ c_1^u f_1^v  m^{u+v}  } \prod_{p \mid \ell q} \Big(L_p(v, \tilde{F} \times \psi)  L_p(u, \tilde{F} \times \bar{\psi})\Big)^{-1}$$
  has continuation to an $\varepsilon$-neighbourhood of $\Re u, \Re v \geq 0$, $\Re(u+v) \geq 1/2$ and is bounded by $(\ell q)^{\varepsilon}$ in this region. (Note that presently $\theta = 0$ since $F = {\tt E}_{0}$.) We estimate the remaining factor on the left hand side of \eqref{bound-cor-eis} trivially by
  $$ (\ell q)^{\varepsilon}  \sum_{\ell_1\ell_2 = \ell} \frac{1}{\ell Q}  \sum_{d_0 \mid M_1 \mid d_0^{\infty}} \sum_{ d_0 M_1M_2 \mid \ell_2Q} \sum_{\delta_1, \delta_2 \mid M_2} \frac{\delta_1 + \delta_2}{M_2}  \ll \frac{1}{\ell Q} (\ell q)^{\varepsilon}$$
for $t = \pm i(s-1)$ and $1/2 -\varepsilon \leq \Re s, \Re w < 1$.  The possible pole at $s = 1/2$ of order at most $\omega(\ell)$ comes from the factor $L(\text{triv}, t, \ell_2)$ at $t = \pm i(1-s)$.

 
\section{The reciprocity formula}

This section is devoted to a proof of Theorem \ref{thm2}. 
Let us initially assume
\begin{equation}\label{assume}
 3/2 < \Re s < 2, \quad 14 < \Re w < 15.
 \end{equation}
In this range we can write the $L$-functions in the definition of $\mathcal{N}_{q, \ell}^{\text{cusp}}(s, w; h)$ as a Dirichlet series, and we have by \eqref{local-RS} and the familiar ${\rm GL}(2)$ Hecke relations
\begin{displaymath}
\begin{split}
& \lambda_f(q) L(s, f \times F) L^{(q)}(w, \bar{f}) \frac{\overline{\Lambda_f}(\ell; w)}{\ell^w}\\
 & = \lambda_f(q) \sum_{n_1, n_2 } \frac{A(n_2, n_1) \chi(n_1) \lambda_f(n_2)}{n_2^s n_1^{2s}} \sum_{(r, q) = 1} \frac{\overline{\lambda_f}(r)}{r^w} \frac{1}{\ell^w} \sum_{\ell_1\ell_2 = \ell} \frac{\mu(\ell_1)\bar{\chi}(\ell_1)\overline{ \lambda_f}(\ell_2)}{\ell_1^w}\\
&=\sum_{n_1, n_2 } \frac{A(n_2, n_1) \chi(n_1) \lambda_f(n_2q)}{n_2^s n_1^{2s}} \sum_{\substack{(r, q) = 1 \\ \ell \mid r}} \frac{\overline{\lambda_f}(r)}{r^w}.
\end{split}
\end{displaymath}
Similarly, since $\lambda_{\mathfrak{a}, t}(n) = \lambda_{\mathfrak{a}, t}(nq)$ by \eqref{lambda-eis} for $\mathfrak{a} \in \{\infty, 0\}$, $t \in \Bbb{R}$ and $n \in \Bbb{N}$, we have  
\begin{displaymath}
\begin{split}
&  L(s, \mathscr{E}_{\mathfrak{a}}(t) \times F) L^{(q)}(w, \overline{ \mathscr{E}_{\mathfrak{a}}(t)}) \frac{\overline{\Lambda_{\mathfrak{a}, t}}(\ell; w)}{\ell^w}=\sum_{n_1, n_2 } \frac{A(n_2, n_1) \chi(n_1) \lambda_{\mathfrak{a}, t}(n_2q)}{n_2^s n_1^{2s}} \sum_{\substack{(r, q) = 1 \\ \ell \mid r}} \frac{\overline{\lambda_{\mathfrak{a}, t}}(r)}{r^w}.
\end{split}
\end{displaymath}
 Hence we can apply the opposite sign Kuznetsov formula \eqref{kuz-neben} to obtain
\begin{displaymath}
  \mathcal{N}_{q, \ell}(s, w; h) =  q\sum_{q \mid c}     \sum_{n_1, n_2} \sum_{\substack{\ell \mid r\\ (r, q) = 1}} \frac{A(n_2, n_1) \chi(n_1)S_{\chi}(-r, n_2q, c)}{n_2^s n_1^{2s} r^wc} \mathscr{K}h\left(\frac{\sqrt{n_2qr}}{c}\right).
\end{displaymath}
We   note that  $S_{\chi}(r, n_2q, c) = 0$ for $(r, q) = 1$, $q^2\mid c$. Indeed, if $c = q^{\alpha} c'$ with $(c', q) = 1$, $\alpha \geq 2$, we have 
$S_{\chi}(r, n_2q, c) = S_{\chi}(r\bar{c}', n_2 q\bar{c}', q^{\alpha}) S(r\bar{q}^{\alpha}, n_2\bar{q}^{\alpha-1}, c')$
with
$$S_{\chi}(r\bar{c}', n_2 q\bar{c}', q^{\alpha})  = \underset{x\, (\text{mod }q^{\alpha-1})}{\left.{\sum} \right.^{\ast}} \sum_{y \, (\text{mod }q)}  \chi(x) e\Big(\frac{(x + q^{\alpha-1} y)r \bar{c}' + n_2 q\bar{c}' \overline{(x + q^{\alpha-1} y)}}{q^{\alpha}} \Big),$$
and the $y$-sum vanishes.  Hence we can write $c = qc'$ with $(q, c') = 1$.  By twisted multiplicativity we get
$$   \mathcal{N}_{q, \ell}(s, w; h)=   \sum_{c, n_1, n_2} \sum_{\substack{\ell \mid r\\ (r, q) = 1}} \frac{A(n_2, n_1) \chi(n_1)S_{\chi}(-r\bar{q}, n_2, c)\tau(\chi) \bar{\chi}( r\bar{c})}{n_2^s n_1^{2s} r^wc} \mathscr{K}h\left(\frac{\sqrt{n_2r}}{q^{1/2}c}\right).$$
Here we also use that $\chi$ is primitive to evaluate the Gau{\ss} sum. At this point we can drop the condition $(r, q) = 1$ which is now automatic. The weight function $\mathscr{K}h$ satisfies the  properties given in Lemma \ref{BLM-lemma}. By Mellin inversion and \eqref{defE} we have
$$ \mathcal{N}_{q, \ell}(s, w; h) = \tau(\chi)\int_{(-1)}\widehat{ \mathscr{K}h}(u) q^{u/2} \mathcal{E}^-_{q, \ell}(s, u, w) \frac{du}{2\pi i} .$$
For $s, w$ as in \eqref{assume}, the triple $(s, u, w)$ is in \eqref{region1}. 
 We shift the $u$-contour to the left (inside the region \eqref{variables}   into the region \eqref{abs1}, namely to  $\Re u = -4/3 - 2\Re s$,  so that we can apply the functional equation of Lemma \ref{func3} and open $\tilde{\mathcal{E}}^{\pm}_{q, \ell}(s, u, w)$ as a Dirichlet series. If $F = {\tt E}_{0}$,  we pick up a pole on the way  whose residue contributes
\begin{equation}\label{Nsw}
\begin{split}
\mathcal{P}_{q, \ell}(s, w; h) &:=    \tau(\chi)  \underset{u = 2-2s}{\text{res}} \widehat{ \mathscr{K}h}(u) q^{u/2} \mathcal{E}^-_{q, \ell }(s, u, w) \\
&= \tau(\chi) \sum_{j=1}^3\frac{1}{(j-1)!} \mathcal{R}_{q, \ell; j}(s, w) \frac{d^{j-1}}{du^{j-1}}  \widehat{ \mathscr{K}h}(u) q^{u/2}\Big|_{u = 2-2s}. 
\end{split}
\end{equation}
It follows from Lemma \ref{BLM-lemma} and Lemma \ref{res} that $\mathcal{P}_{q, \ell}(s, w; h)$ has meromorphic continuation to an $\varepsilon$-neighbourhood of $\Re s, \Re w \geq 1/2 $. In an $\varepsilon$-neighbourhood of the region  $1/2 \leq \Re s \leq \Re w$ the only pole can occur at $s+w= 1$,  and the bound 
\begin{equation}\label{bound-P}
(s+w-1)^4 \mathcal{P}_{q, \ell}(s, w; h)\ll_{s, w}   q^{3/2 - \Re s} T^{3-2\Re s}   (\ell^{-\Re(s+w)} + \ell^{-\Re(1+w-s)}) (q\ell T)^{\varepsilon} \ll qT^2\ell^{-1} (q\ell T)^{\varepsilon}
\end{equation}
holds for any $\varepsilon > 0$. 

The shifted integral equals
\begin{displaymath}
\begin{split}
 \mathcal{N}^{\ast}_{q, \ell}(s, w; h) =   \tau(\chi) \int_{(-4/3 - 2\Re s)}  & \widehat{ \mathscr{K}h}(u) q^{u/2} 
\sum_{\pm} \mathcal{G}_{-\mu}^{\pm}(1-s-u/2) \\
& \sum_{(c, d) = 1} \sum_{\substack{(m, q) = 1\\ \ell \mid dm}} \sum_n  \frac{A(m, n) \bar{\chi}(d)\chi(c) e(\pm \bar{d} n q/c)}{m^{s+w} n^{1-s-u/2} c^{-1+3s+u/2} d^{w+u/2}}
\frac{du}{2\pi i}. 
\end{split}
\end{displaymath}
At this point we insert artificially a factor
$$1 = e\left(\mp \frac{nq}{cd}\right) e\left(\pm \frac{nq}{cd}\right) = e\left(\mp \frac{nq}{cd}\right) \int_{\mathcal{C}}  G^{\pm}(\xi)\left(\frac{nq}{cd}\right)^{-\xi} \frac{d\xi}{2\pi i},$$
with $G^{\pm}$ as in \eqref{defG} and  $\mathcal{C} = \mathcal{C}(1/10, -3/5)$, where $\mathcal{C}(x, y)$ is the curved contour 
 \begin{equation*}
   \textstyle (y -i\infty, y - i] \cup [y - i, x] \cup [x, y + i] \cup [y + i, y + i\infty). 
 \end{equation*} 
   By additive reciprocity, we obtain
 \begin{equation}\label{rec}
\begin{split}
\mathcal{N}^{\ast}_{q, \ell}(s, w; h) =  \tau(\chi) & \int_{(-4/3 - 2\Re s)}   \int_{\mathcal{C}} \widehat{ \mathscr{K}h}(u) q^{u/2+\xi}\sum_{\pm} \mathcal{G}_{-\mu}^{\pm}(1-s-u/2)  G^{\pm}(\xi) \\
&  \sum_{\substack{(m, q) = 1\\ \ell \mid dm}} \sum_n  \frac{A(m, n) \bar{\chi}(d) \Theta_q( \pm n, d, -1+3s+u/2 - \xi)}{m^{s+w} n^{1-s-u/2+\xi}   d^{w+u/2-\xi}}
 \frac{d\xi}{2\pi i} \frac{du}{2\pi i}
\end{split}
\end{equation}
with $\Theta_q(n, d; v)$ as in \eqref{defTheta}.  
 Since $\Re(1-s-u/2 + \xi) = \Re(5/3 - \xi) \geq 16/15 > 1$ as well as $\Re(-1+3s+u/2 - \xi) = \Re(-5/3 + 2s -\xi) > 37/30 > 1$ for $\Re u = -4/3 - 2\Re s$, $\Re s > 3/2$ and $\xi \in \mathcal{C}(1/10, -3/5)$, the expression \eqref{rec} (even after opening $\Theta_q$) is absolutely convergent.  
We interchange the two integrals, slightly curve the inner $u$-integral so that $u' = u- 2\xi$ becomes a straight line $\Re u' =-4/3 - 2\Re s$ getting 
 \begin{equation*}
\begin{split}
\mathcal{N}^{\ast}_{q, \ell}(s, w; h) = \tau(\chi)  \int_{\mathcal{C}} &  \int_{(-4/3 - 2\Re s)}    \widehat{ \mathscr{K}h}(u'+2\xi) q^{u'/2} 
\sum_{\pm} \mathcal{G}_{-\mu}^{\pm}(1-s-u'/2-\xi)  G^{\pm}(\xi) \\
&  \sum_{\substack{(m, q) = 1\\ \ell \mid dm}} \sum_n  \frac{A(m, n) \bar{\chi}(d) \Theta_q( \pm n, d, -1+3s+u'/2)}{m^{s+w} n^{1-s-u'/2}   d^{w+u'/2}}
 \frac{du'}{2\pi i}  \frac{d\xi}{2\pi i},
\end{split}
\end{equation*}
 which equals
$$ \tau(\chi)  \int_{\mathcal{C}}  \int_{(-4/3 - 2\Re s)}    \widehat{ \mathscr{K}h}(u'+2\xi)
q^{u'/2} 
\sum_{\pm} \mathcal{G}^{\pm}_{-\mu}(1-s-u'/2-\xi)  G^{\pm}(\xi)  \mathcal{A}^{\pm}(s, u', w)
 \frac{du'}{2\pi i}  \frac{d\xi}{2\pi i}
$$
with the  notation \eqref{defAB}. The triple $(s, u', w)$ is now in \eqref{abs1}. Now we continue to shift the $u'$-contour further to the left in order to move into the region \eqref{b}. This is not directly possible because the $\xi$-integral will no longer converge. In order to get better convergence properties, we bend the $\xi$-contour $\mathcal{C} = \mathcal{C}(1/10, -3/5)$ more and replace it with $\mathcal{C}(1/10, -3/2)$. This does not cross any poles. Now we can shift the $u$-contour to the left  to $\Re u' = 7/4 - 6\Re s$.  Condition \eqref{assume} ensures that we end up in \eqref{b}, and by \eqref{414a}, Stirling's formula for $G^{\pm}(\xi)$ and Lemma \ref{BLM-lemma} the double integral is absolutely convergent.  
This  does not cross any poles since the argument of the gamma factors moves to the right.  We apply the functional equation \eqref{func-new} 
getting
\begin{equation*}
\begin{split}
\mathcal{N}^{\ast}_{q, \ell}(s, w; h) = \frac{\tau(\chi)^2}{q^{3s-1}}&   \int_{\mathcal{C}}  \int_{(7/4 - 6\Re s)}       \widehat{ \mathscr{K}h}(u'+2\xi) 
\sum_{\pm} \mathcal{G}_{-\mu}^{\pm}(1-s-u'/2-\xi)  G^{\pm}(\xi)   \\
&\sum_{\epsilon \in \{\pm 1\}} G^{\mp \epsilon}(2 - 3s - u'/2) \mathcal{B}^{\pm \epsilon}(s, u', w)
 \frac{du'}{2\pi i}  \frac{d\xi}{2\pi i}.
\end{split}
\end{equation*}
In this region we can open the $\mathcal{B}$-function, i.e.\ 
 $$\mathcal{B}^{\pm \epsilon}(s, u', w) =  \sum_{\substack{(md, q) = 1\\ \ell \mid dm}} \sum_n \sum_c \frac{  \bar{\chi}(c) S(\pm \epsilon n,   c, d)A(m, n)   }{c^{2-3s-u'/2} m^{s+w} n^{1-s-u'/2}   d^{3s+u'+w-1}}.$$
Note that the double integral is  absolutely convergent by \eqref{coro-bound} and the choice of the $\xi$-contour. 
We introduce a new set of variables
\begin{equation}\label{newvar}
3s+w+u'-1 = 1-u'', \quad (s+w)/2 = s', \quad 1+w/2 - 3s/2 = w'
\end{equation}
getting
\begin{displaymath}
\begin{split}
 \frac{   \tau(\chi)^2}{q^{(1+ 3s'-3w')/2}}&   \int_{\mathcal{C}}  \int_{(9/4 - 2\Re w')}      \widehat{ \mathscr{K}h}(1+w'-3s' - u'' +2\xi) 
\sum_{\pm} \mathcal{G}_{-\mu}^{\pm}(s' + u''/2 - \xi)  G^{\pm}(\xi)   \\
&\sum_{\epsilon \in \{\pm 1\}} G^{\mp \epsilon}(w' + u''/2)   \sum_{\substack{(md, q) = 1\\ \ell \mid dm}} \sum_n \sum_c \frac{\bar{\chi}(c)  S(\pm \epsilon n,   c, d)A(m, n)   }{c^{w' + u''/2} m^{2s'} n^{s' + u''/2}   d^{1-u''}}
 \frac{du''}{2\pi i} \frac{d\xi}{2\pi i} .
\end{split}
\end{displaymath}
Changing integrals once again we finally arrive at
\begin{equation}\label{finally}
\begin{split}
  \mathcal{N}^{\ast}_{q, \ell}(s, w; h) =  \frac{   \tau(\chi)^2}{q^{(1+ 3s'-3w')/2}}&  \sum_{\pm} \sum_{\epsilon \in \{\pm 1\}}   \int_{(9/4 - 2\Re w')}    \mathcal{H}^{\pm, \epsilon}_{s', w'}(u'') \\
 &     \sum_{\substack{(md, q) = 1\\ \ell \mid dm}} \sum_n \sum_c \frac{  \chi(\bar{c}) S(\pm \epsilon n,   c, d)A(m, n)   }{c^{w' } m^{2s'} n^{s' }   d}\left(\frac{\sqrt{nc}}{d}\right)^{-u''}
 \frac{du''}{2\pi i} .
\end{split}
\end{equation}
where
\begin{equation}\label{defHtrafo}
\mathcal{H}^{\pm, \epsilon}_{{\tt s}, {\tt w}}({\tt u})  =  \int_{\mathcal{C}}     \widehat{ \mathscr{K}h}(1+{\tt w}-3{\tt s} - {\tt u} +2\xi) 
  \mathcal{G}_{-\mu}^{\pm}({\tt s} + {\tt u}/2 - \xi)  G^{\pm}(\xi)G^{\mp \epsilon}({\tt w} + {\tt u}/2) \frac{d\xi}{2\pi i} 
  \end{equation}
with $\mathcal{C} = \mathcal{C} (1/10, -3/2)$. 

\begin{lemma}\label{H-transform} Let $T \geq 1$, $B > 0$ and let $h$ be a $T$-admissible function. If $A$ in \eqref{admi} is sufficiently large in terms of $B$, then the function $\mathcal{H}^{\pm, \epsilon}_{{\tt s}, {\tt w}}({\tt u}) $  is meromorphic in $({\tt s}, {\tt w}, {\tt u})$ in the region 
\begin{equation}\label{region}
 \Re({\tt s} + {\tt u}/2)> -1+\theta , \quad \Re({\tt w} + {\tt u}/2) > -1, \quad |\Re {\tt s}|, |\Re {\tt u}|, |\Re {\tt w}| \leq B
 \end{equation}
with poles at most at ${\tt s} + {\tt u}/2-\mu_j= 0$, ${\tt w} + {\tt u}/2 = 0$,  
and satisfies
\begin{equation}\label{h-trafo-bound}
\mathcal{H}^{\pm, \epsilon}_{{\tt s}, {\tt w}}({\tt u})\Big[({\tt w} + {\tt u}/2) \prod_{j=1}^3({\tt s} + {\tt u}/2 - \mu_j) \Big] \ll_{{\tt s}, {\tt w}} T^{2 + \Re({\tt w} - 3{\tt s} -   {\tt u}) }(1+|{\tt u}|)^{-B}.
\end{equation}
\end{lemma}

\textbf{Proof.} We recall that $ \widehat{ \mathscr{K}h}$ is holomorphic in a wide strip and rapidly decaying by Lemma \ref{BLM-lemma}.   We replace the contour in \eqref{defHtrafo} once again by $\mathcal{C}(-1+\varepsilon, -B_0)$ for some very large $B_0$ and some sufficiently small $\varepsilon > 0$. We pick up a pole at $\xi = 0$ with residue $$\widehat{ \mathscr{K}h}(1+{\tt w}-3{\tt s} - {\tt u} )  \mathcal{G}_{-\mu}^{\pm}({\tt s} + {\tt u}/2  )   G^{\mp \epsilon}({\tt w} + {\tt u}/2) $$
which is holomorphic in \eqref{region} except for poles at ${\tt s} + {\tt u}/2 - \mu_j = 0$, ${\tt w} + {\tt u}/2 = 0$ and it satisfies the bound \eqref{h-trafo-bound} if $A$ is sufficiently large. The remaining integral is holomorphic in \eqref{region} except for a pole at ${\tt w} + {\tt u}/2 = 0$ and can be estimated   trivially using \eqref{coro-bound} by
$$ \ll_{{\tt s}, {\tt w}} T^{1+\varepsilon + \Re({\tt w} - 3{\tt s} -   {\tt u})} (1 + |{\tt u}|) ^{\Re ({\tt w} +{\tt u}/2 - 1/2)} \int_{\mathcal{C}(-1+\varepsilon, -B_0)} (1 +  |\Im {\tt u} + 2\xi|)^{-B_0}(1 + |\xi|)^{\Re \xi -1/2} |d\xi|.$$
 Taking $B_0$ sufficiently large (say $B_0 = 100B$), the desired bound follows.  \\
 
 We return to \eqref{finally} and open the condition $\ell \mid dm$. In this way we see that the $m, d, n, c$-sum equals 
$$\sum_{\ell_1\ell_2 = \ell} \frac{1}{\ell_1^{2s'}}\sum_{ \substack{(d, q) = 1\\ \ell_2\mid d}} \sum_{(m, \ell_2 q) = 1}\sum_n \sum_c \frac{ \bar{\chi}(c)  S(\pm \epsilon n,   c, d)A(\ell_1m, n)   }{c^{w' } m^{2s'} n^{s' }   d} .$$
Now taking the inverse Mellin transform (with notation as in Lemma \ref{lem3}) and applying the Kuznetsov formula \eqref{kuz2}, we obtain
\begin{displaymath}
\begin{split}
  \mathcal{N}^{\ast}_{q, \ell}(s, w; h) =  \frac{   \tau(\chi)^2}{q^{(1+ 3s'-3w')/2}}  \sum_{\pm} \sum_{\epsilon \in \{\pm 1\}}& \sum_{\ell_1\ell_2 = \ell} \frac{1}{\ell_1^{2s'}}  \sum_{(m, \ell_2 q) = 1}\sum_n \sum_c \frac{ \bar{\chi}(c)  A(\ell_1m, n)   }{c^{w' } m^{2s'} n^{s' }   d} \\
  &\left( \mathcal{A}_{\ell_2}(\pm \epsilon n, c, \mathscr{L}_{\pm \epsilon}\widecheck{ \mathcal{H}^{\pm, \epsilon}_{s', w'}}) -  \mathcal{A}_{\ell_2q}(\pm \epsilon n, c,\mathscr{L}_{\pm \epsilon}\widecheck{ \mathcal{H}^{\pm, \epsilon}_{s', w'}} )\right).
\end{split}
\end{displaymath}
The analytic properties of $\mathscr{L}_{\pm \epsilon}\widecheck{ \mathcal{H}^{\pm, \epsilon}_{s', w'}} $ were studied in Lemma \ref{lem3}. Note that $ \mathcal{H}^{\pm, \epsilon}_{s', w'}$ satisfies the assumptions of that lemma by Lemma \ref{H-transform}  for some arbitrarily  large $B$, $$c = 2 + \Re (w '-3s'), \quad \rho = 1 + \Re s' - \theta,$$ and poles at most at $u = -2w$, $u = 2\mu_j - 2s'$, $j = 1, 2, 3$. In particular, $\mathscr{L}_{\pm \epsilon}\widecheck{ \mathcal{H}^{\pm, \epsilon}_{s', w'}} $ has meromorphic continuation to $|\Im t| < \min(1 + \Re s' - \theta, 1 + \Re w')$ with poles at most at 
\begin{equation}\label{poles}
\pm iw', \quad \pm i(\mu_j - s'), \quad j = 1, 2, 3,
\end{equation} and 
the bounds \eqref{L1} for the non-exceptional and Eisenstein spectrum, \eqref{L3} for the holomorphic spectrum and \eqref{L2} for the (possible) exceptional spectrum (as long as $ \vartheta_0 + \theta \leq 7/64 + 5/14 < 1/2 \leq \Re s'$) are applicable, so that the  above spectral sum
 is absolutely (and in fact rapidly) convergent.  Inserting \eqref{AtoM1} and \eqref{AtoM2}, we obtain finally 
$$    \mathcal{N}_{q, \ell}(s, w; h) =  \mathcal{P}_{q, \ell}(s, w; h) + \frac{   \tau(\chi)^2}{q^{(1+ 3s'-3w')/2}}  \sum_{\pm} \sum_{\epsilon \in \{\pm 1\}} \left(\mathcal{M}^{\pm \epsilon}_{\ell}(s', w';\mathcal{T}_{s', w'}^{\pm, \epsilon}h) - \mathcal{M}^{\pm \epsilon}_{\ell q}(s', w'; \mathcal{T}_{s', w'}^{\pm, \epsilon}h) \right)$$
under the assumption   \eqref{assume} where
\begin{equation}\label{defT}
\mathcal{T}_{s', w'}^{\pm, \epsilon}h = \mathscr{L}_{\pm \epsilon}\widecheck{ \mathcal{H}^{\pm, \epsilon}_{s', w'}}
\end{equation}
with $\widecheck{ \mathcal{H}^{\pm, \epsilon}_{s', w'}}$ as in \eqref{defHtrafo}. 

 Analytic continuation to an $\varepsilon$-neighbourhood of $1/2 \leq \Re s \leq \Re w < 1$ (these inequalities imply also $1/2 \leq \Re s', \Re w' < 1$) follows now as in \cite[Section 10]{BK}.  Lemma \ref{lem3} ensures that the transform $\mathcal{T}_{s', w'}^{\pm, \epsilon}$ can be continued and also provides the bound \eqref{weakly} to show that $\mathcal{T}_{s', w'}^{\pm, \epsilon}$ is weakly $T$-admissible for $\Re s' = \Re w' = 1/2$ (equivalently $\Re s = \Re w = 1/2$).  The continuation of $\mathcal{P}_{q, \ell}(s, w; h)$ (if $F = {\tt E}_{0}$) with poles at most at $s+w= 1$ was discussed at \eqref{Nsw}. 

The analytic continuation of $ \mathcal{N}_{q, \ell}(s, w; h) $  to an $\varepsilon$-neighbourhood of $1/2 \leq \Re s, \Re w < 1$  incurs by Lemma \ref{ana-cont}a the  additional (holomorphic) polar term $\mathcal{R}_{q, \ell}(s, w; h)$  if $F = {\tt E}_{0}$. 

Finally the analytic continuation of $\mathcal{M}^{\pm \epsilon}_{\ell Q}(s', w'; \mathcal{T}_{s', w'}^{\pm, \epsilon}h)$ for $Q \in \{1, q\}$ incurs by   Lemma \ref{ana-cont}b the  additional polar term   $\widetilde{\mathcal{R}}_{\ell Q}(s', w'; \mathcal{T}_{s', w'}^{\pm, \epsilon}h)$  if $F = {\tt E}_{0}$. This lemma is applicable because $(\ell, q) = 1$, so that $q^2 \nmid \ell Q$, and the decay/regularity properties of the weight function are ensured by Lemma \ref{lem3}. By \eqref{poles} and Lemma \ref{ana-cont}b, these terms can have poles at most at $s' = 1/2$ and $s'=(1- w')$ in $1/2 \leq \Re s', \Re w' < 1$.  Combining \eqref{L2} with $|\Im t| = |\Re s' - 1|$, $c = 2 + \Re (w '-3s')$, \eqref{bound-cor-eis} for $\theta = 0$ (since $F = {\tt E}_{\mu}$) and the convexity bound with \eqref{cs} we obtain
\begin{equation*}
\begin{split}
(s' - 1/2)^{\omega(\ell) + 3} & (s' + w' - 1)\frac{   \tau(\chi)^2}{q^{(1+ 3s'-3w')/2}}  \widetilde{\mathcal{R}}_{\ell Q}(s', w'; \mathcal{T}_{s', w'}^{\pm, \epsilon}h)\\
& \ll_{s', w'} (q\ell T)^{\varepsilon} (\ell Q)^{-1} T^{4 + \Re (w'- 5 s')}   q^{c(s'+w'-1) + c(1+w'-s') + 3(s'-w') + 1/2}.
\end{split}
\end{equation*}
In particular, if $1/2 - \varepsilon <  \Re s,   \Re w < 1/2+\varepsilon$, so that $1/2 - \varepsilon <  \Re s',   \Re w' < 1/2+\varepsilon$, then by $c(0) = 1/2$, $c(1) = 0$ we obtain
\begin{equation}\label{bound-tildeR}
(s' - 1/2)^{\omega(\ell) + 3}  (s' + w' - 1) \frac{   \tau(\chi)^2}{q^{(1+ 3s'-3w')/2}} \widetilde{\mathcal{R}}_{\ell Q}(s', w'; \mathcal{T}_{s', w'}^{\pm, \epsilon}h) \ll_{ s', w'} (q\ell T)^{\varepsilon} \frac{T^2   q}{\ell Q}.
\end{equation}
This gives the desired reciprocity formula of Theorem \ref{thm2} with
\begin{equation}\label{defGmain}
\begin{split}
\mathcal{G}_{q, \ell}(s, w; h) = & \mathcal{P}_{q, \ell}(s, w; h)   -\mathcal{R}_{q, \ell}(s, w; h) \\
&+ \frac{   \tau(\chi)^2}{q^{(1+ 3s'-3w')/2}}  \sum_{\pm} \sum_{\epsilon \in \{\pm 1\} }\left(\widetilde{\mathcal{R}}_{\ell}(s, w; \mathcal{T}_{s', w'}^{\pm, \epsilon}h) - \widetilde{\mathcal{R}}_{\ell q}(s, w; \mathcal{T}_{s', w'}^{\pm, \epsilon}h)\right) .
\end{split}
\end{equation}
 While the individual terms may have poles as described above, the expression $\mathcal{G}_{q, \ell}(s, w; h) $ must be holomorphic, since the rest of the formula in Theorem \ref{thm2} is holomorphic. Note that the zeros imposed on  $h$ in Theorem \ref{thm2} imply that $\mathcal{R}_{q, \ell}(s, w; h) = 0$, cf.\ \eqref{eis-expl}.  The bound  \eqref{bound-G} follows now from \eqref{bound-P} and \eqref{bound-tildeR} outside of $s' = 1/2$, $s' + w' = 1$ and $s+w=1$ and also in a neighbourhood of these lines by an application of a two-dimensional Cauchy integral formula at the cost of a factor   $\ll \varepsilon^{-\omega(\ell) } \ll \ell^{\varepsilon}$.  This completes the proof of Theorem \ref{thm2}.

 \section{Proof of Theorem \ref{thm3}}
 
Before we start with the proof, we recall that a standard application of the large sieve (cf.\ e.g.\ \cite[Theorm 7.35]{IK}) shows that
\begin{equation}\label{largesieve}
\sum_{f \in \mathcal{B}^{\ast}(N, \text{triv})}  L(1/2, f)^4 (1 + |t_f|)^{-10}+ 
    \sum_{f \in   \mathcal{B}_{\text{hol}}^{\ast}(N, \text{triv})}  L(1/2, f)^4  k_f^{-10}  \ll N^{1+\varepsilon}
    \end{equation}
    for any $N \in \Bbb{N}$. Moreover, by another standard application of the large sieve (\cite[Theorem 7.34]{IK}) we have 
   \begin{equation}\label{largesieve2}
   \sum_{c_{\psi}^2 \mid N} \int_{-\infty}^{\infty} |L(1/2 + it, \psi)|^8  (1 + |t|)^{-10} dt \leq   \sum_{c_{\psi}\leq N^{1/2} } \int_{-\infty}^{\infty} |L(1/2 + it, \psi)|^8   (1 + |t|)^{-10} dt  \ll N^{1+\varepsilon}.
      \end{equation} 
    The bound
    \begin{equation}\label{zeta4}
    \int_{-T}^T |\zeta(1/2 + it)|^4 dt \ll T^{1+\varepsilon}
    \end{equation}  
is classical.    We will also need the following mean value result
\begin{equation}\label{mean-value}
  \sum_{d \mid 3[N, q]} \sum_{f \in \mathcal{B}^{\ast}(d, \text{triv})} \frac{|L(w, f \times \chi)|^2}{(1 +|t_f|)^r} \ll_w [N, q]^{1+\varepsilon} 
\end{equation}
for $\Re w = 1/2 + \varepsilon$ and $r > 6$ where $[N, q]$ denotes the least common multiple of $N$ and $q$. This follows from \cite[(3.9)]{BH} in the special case $\ell = 1$ in combination with \cite[(2.7), (3.1), (3.3), (3.7)]{BH}. By the functional equation and the Phragm\'en-Lindel\"of principle, the bound remains true for $\Re w = 1/2$, upon changing the value of $\varepsilon$ on the right hand side of \eqref{mean-value}.

    We recall the standard bounds
\begin{equation}\label{hl}
\begin{split}
(N(1 + |t_f|))^{-\varepsilon} & \ll L(1, \text{Ad}^2f) \ll (N(1 + |t_f|))^{\varepsilon}, \quad f \in \mathcal{B}^{\ast}(N),\\
 (N k_f)^{-\varepsilon}  &\ll L(1, \text{Ad}^2f) \ll (N k_f)^{\varepsilon}, \quad\quad\quad\,\,\,\,\, f \in \mathcal{B}_{\text{hol}}^{\ast}(N),\\
 (c_{\psi}(1+|t|))^{-\varepsilon} & \ll |L(1 + 2 it, \psi)|.
 \end{split}
\end{equation}
 Next we recall the subconvexity bound \cite[Theorem 2]{BH} 
\begin{equation}\label{BH}
L(s, f\times  \chi) \ll_{s} (1+|t_f|)^4 (q^{3/8} N^{1/4} + q^{1/4} N^{1/2})(Nq)^{\varepsilon}
\end{equation}
for $f \in \mathcal{B}^{\ast}(N, \text{triv})$, $\Re s = 1/2$  and $(N, q) = 1$, and analogously for $f \in \mathcal{B}_{\text{hol}}^{\ast}(N, \text{triv})$. We also have the hybrid subconvexity bound \cite{HB}
\begin{equation}\label{Heath-Brown}
L(1/2 + it, \psi) \ll ((1+|t|)c_{\psi})^{3/16+\varepsilon}
\end{equation}
for $t \in \Bbb{R}$ and a character $\psi$ of conductor $c_{\psi}$. We start with the following key lemma.
\begin{lemma} Let $\Re s = \Re w = 1/2$, $F = {\tt E}_{0}$, $(\ell, q) = 1$ and $\mathfrak{h} = (h, h^{\text{{\rm hol}}})$ weakly $T$-admissible as in \eqref{weakly}. Then
\begin{equation}\label{M1}
\mathcal{M}^{\pm}_{\ell}(s, w; \mathfrak{h}) \ll_{s, w} (q^{3/8} \ell^{1/4} + q^{1/4} \ell^{1/2})T^{1+2\vartheta_0} (\ell qT)^{\varepsilon}.
\end{equation}
and
\begin{equation}\label{M2}
\mathcal{M}^{\pm}_{\ell q}(s, w; \mathfrak{h}) \ll_{s, w} T^{1+2\vartheta_0} (\ell q)^{1/4} (\ell q T)^{\varepsilon}
\end{equation}
for any $\varepsilon > 0$. 
\end{lemma}

\textbf{Proof.} We bound 
   $\mathcal{M}^{\text{Maa{\ss}}, \pm }_{\ell}(s, w; \mathfrak{h})  $ trivially by
$$\sum_{\ell_0 \mid \ell} \sum_{f \in \mathcal{B}^{\ast}(\ell_0, \text{triv})} \frac{|L(w, f\times  \bar{\chi})    L(s , f)^3 | }{L(1, \text{Ad}^2 f)} \frac{(\ell q)^{\varepsilon}}{\ell} T^{1+2\vartheta_0+\varepsilon} (1 + |t_f|)^{-20}$$
using \eqref{corr-bound-final} with $\theta = 0$. We use \eqref{hl} for $L(1, \text{Ad}^2f)$, \eqref{BH} for $L(w, f\times \bar{\chi})$ and \eqref{largesieve} with $N = \ell_0$  for the rest to obtain
$$\mathcal{M}^{\text{Maa{\ss}}, \pm }_{\ell}(s, w; \mathfrak{h})  \ll _{s, w} (q^{3/8} \ell^{1/4} + q^{1/4} \ell^{1/2})T^{1+2\vartheta_0} (\ell qT)^{\varepsilon}. $$
 The same bound holds (even without $T^{2\vartheta_0+\varepsilon}$) for the corresponding holomorphic term. We estimate the Eisenstein term trivially by (using \eqref{corr-bound-final} with $\theta = 0$)
$$\sum_{\psi : c_{\psi}^2 \mid \ell} \int_{\Bbb{R}} \frac{ |L(w +it, \bar{\chi}\psi)L(w -it, \overline{\chi\psi})   L(s +it , \psi)^3 L(s - it , \bar{\psi})^3| }{|L(1 + 2it, \psi^2)|^2(1+|t|)^{20}  } \frac{(\ell q)^{\varepsilon}T}{\ell} dt.$$
Here we apply \eqref{Heath-Brown} for   $L(w +it, \bar{\chi}\psi)L(w -it, \overline{\chi\psi}) $ (recall that 
 $c_{\psi} \leq \ell^{1/2}$) and \eqref{largesieve2} for the rest along with \eqref{hl} to confirm   $$
\mathcal{M}^{\text{Eis}}_{\ell}(s, w; \mathfrak{h}) \ll_{s, w} T q^{3/8} \ell^{3/16} (\ell qT)^{\varepsilon}.$$
Since $(\ell, q) = 1$, the condition $c_{\psi}^2 \mid \ell q$ is equivalent to $c_{\psi}^2 \mid \ell$, so that the same argument gives
$$
\mathcal{M}^{\text{Eis}}_{\ell q}(s, w; \mathfrak{h}) \ll_{s, w} T q^{-5/8} \ell^{3/16} (\ell qT)^{\varepsilon}.$$
This completes the proof of \eqref{M1} and part of \eqref{M2}, and it remains to estimate the cuspidal contribution of  $\mathcal{M}^{  \pm }_{\ell q}(s, w; \mathfrak{h}) $. Here we use the Cauchy-Schwarz inequality to obtain
\begin{displaymath}
\begin{split}
\mathcal{M}^{ \text{Maa{\ss}}, \pm }_{\ell q}(s, w; \mathfrak{h}) \ll_{s, w}& \frac{T^{1+2\vartheta_0}}{\ell q} (\ell q T)^{\varepsilon}  \left(\sum_{N \mid \ell q} \sum_{f \in \mathcal{B}^{\ast}(N, \text{triv})} \frac{|L(w, f \times \bar{\chi})|^2}{(1 + |t_f|)^{15}}\right)^{1/2}\\
&  \left(\sum_{N \mid \ell q}\sum_{f \in \mathcal{B}^{\ast}(N, \text{triv})} \frac{|L(s , f)|^4}{(1 + |t_f|)^{15}}\right)^{1/2} \max_{\substack{ f \in \mathcal{B}^{\ast}(N, \text{triv}) \\ N \mid \ell q}}   \frac{|L(s , f)|}{(1 + |t_f|)^{4}}.
\end{split}
\end{displaymath}
Now we apply \eqref{mean-value}, \eqref{largesieve} and the convexity bound to obtain
$$\mathcal{M}^{ \text{Maa{\ss}}, \pm }_{\ell q}(s, w; \mathfrak{h}) \ll_{s, w} \frac{T^{1+2\vartheta_0}}{\ell q} (\ell q T)^{\varepsilon}  (\ell q)^{1/2}  (\ell q)^{1/2} (\ell q)^{1/4} \ll T^{1+2\vartheta_0} (\ell q)^{1/4} (\ell q T)^{\varepsilon} ,$$
as desired. The same argument works for the holomorphic part. This completes the proof. \\

For fixed $\tau \in \Bbb{R}$ we choose   $s = \frac{1}{2}+ i\tau$, $w = \frac{1}{2} - i\tau$.  Let $h^{\ast}$ be any $T$-admissible function. Let $(\ell, q) = 1$ and $F = {\tt E}_0$. Then $$\mathcal{N}^{\text{cusp}}_{q, \ell}(s, w; h^{\ast}) = \sum_{f\in \mathcal{B}(q, \chi)} \epsilon_f\lambda_f(q) \frac{L(1/2 + i\tau, f)^3   L^{(q)}(1/2 - i\tau, \bar{f})}{L(1, \text{Ad}^2 f)} \frac{\overline{\Lambda_f}(\ell; 1/2 - i\tau)}{\ell^{1/2 - i\tau}} h^{\ast}(t_f).$$
The functional equation \cite[Proposition 8.1 with $k=0$]{DFI}  along with \cite[Theorem 6.29]{Iw} (which is a purely formal computation to calculate the Fricke eigenvalue and holds also for Maa{\ss} forms) states
$${\tt G}_{\epsilon_f}(1/2 +i\tau, t_f)  L(1/2 + i\tau, f)  = q^{-i\tau -1/2} \epsilon_f\overline{ \lambda_f(q) \tau(\bar{\chi})} {\tt G}_{\epsilon_f}(1/2 - i\tau, t_f) L(1/2 - i\tau, \bar{f})$$
 where
 $${\tt G}_{\epsilon}(s, t) =  \pi^{-s}\Gamma\left(\frac{1}{2}\left(s+\frac{1-\epsilon}{2}-it\right)\right)\Gamma\left(\frac{1}{2}\left(s+\frac{1-\epsilon}{2}+it\right)\right).$$
Thus, for our choice of $s, w$ and $F$ the term  $\mathcal{N}^{\text{cusp}}_{q, \ell}(s, w; h^{\ast})$ equals
$$ \frac{\overline{\tau(\bar{\chi})}}{q^{1/2 + i\tau}}\sum_{f\in \mathcal{B}(q, \chi)}   \frac{L(1/2 + i\tau, f)^2 L(1/2 - i\tau, f) L^{(q)}(1/2 - i\tau, \bar{f})}{L(1, \text{Ad}^2 f)}  \frac{\overline{\Lambda_f}(\ell; 1/2 - i\tau)}{\ell^{1/2 - i\tau}} \frac{{\tt G}_{\epsilon_f}(1/2 - i\tau, t_f) }{ {\tt G}_{\epsilon_f}(1/2 +i\tau, t_f)} h^{\ast}(t_f). $$
By M\"obius inversion we have
$$\lambda_f(\ell) = \sum_{ab = \ell} \Lambda_f(a; w) \frac{\chi(b)}{b^{w}},$$
so that (using also $\overline{\lambda_f(\ell)} = \bar{\chi}(\ell) \lambda_f(\ell)$ for $(\ell, q) = 1$)
\begin{displaymath}
\begin{split}
& \frac{q^{1/2 + i\tau} \ell^{1/2 - i\tau}  \chi(\ell)}{\overline{\tau(\bar{\chi})}}  \sum_{ab = \ell} \frac{\bar{\chi}(b)}{b^{1/2+i\tau}}\mathcal{N}^{\text{cusp}}_{q, a}(s, w; h^{\ast}) \\
 & = \sum_{f\in \mathcal{B}(q, \chi)}   \frac{L(1/2 + i\tau, f)^2 L(1/2 - i\tau, \bar{f})L^{(q)}(1/2 - i\tau, \bar{f}) }{L(1, \text{Ad}^2 f)}  \lambda_f(\ell) \frac{{\tt G}_{\epsilon_f}(1/2 - i\tau, t_f) }{ {\tt G}_{\epsilon_f}(1/2 +i\tau, t_f)} h^{\ast}(t_f).
\end{split}
\end{displaymath}

On the other hand, if $h^{\ast}$ is a $T$-admissible function, 
then we can bound $\mathcal{N}_{q, \ell}^{\text{Eis}}(s, w, h^{\ast})$ for  our choice of $s, w$ and $F$ trivially by inserting \eqref{eis-expl} and  \eqref{Heath-Brown}, so that 
$$\mathcal{N}_{q, \ell}^{\text{Eis}}(s, w, h^{\ast}) \ll  (q\ell T)^{\varepsilon} \frac{(qT)^{3/4 }}{\ell^{1/2}}   \int_{\Bbb{R}} e^{-(t/T)^2} \frac{|\zeta(1/2 \pm  i t + i\tau)|^4   }{|L(1 + 2it, \chi)|^2} dt \ll \frac{T^{7/4} q^{3/4} }{\ell^{1/2}} (q\ell T)^{\varepsilon} $$
by \eqref{zeta4} and \eqref{hl}. (Better estimates could be obtained easily, but the present bound suffices.)

Inserting this in the reciprocity formula of Theorem \ref{thm2}, we obtain under the additional assumption that $h^{\ast}$ has a triple zero at $\pm i(1/2 \pm i\tau)$   that
\begin{displaymath}
\begin{split}
& \sum_{f\in \mathcal{B}(q, \chi)}    \frac{L(1/2 + i\tau, f)^2 L(1/2 - i\tau, \bar{f})L^{(q)}(1/2 - i\tau, \bar{f})}{L(1, \text{Ad}^2 f)}    \lambda_f(\ell) \frac{{\tt G}_{\epsilon_f}(1/2 - i\tau, t_f) }{ {\tt G}_{\epsilon_f}(1/2 +i\tau, t_f)} h^{\ast}(t_f)\\
 &\ll (\ell q T)^{\varepsilon} \left( \frac{qT^2}{\ell^{1/2}} + q^{3/4} T^{7/4} + (\ell q)^{1/2}\sum_{ab = \ell} \frac{1}{b^{1/2}} \sum_{\pm} \sum_{\epsilon \in \{\pm 1\}} \sum_{Q \in \{1, q\}} \Big|\mathcal{M}^{\pm \epsilon}_{aQ}(\textstyle\frac{1}{2} , \frac{1}{2}-2 i \tau;\mathcal{T}_{\frac{1}{2} , \frac{1}{2} - 2i\tau}^{\pm, \epsilon}h)\Big|   \right)\\
 & \ll  (\ell q T)^{\varepsilon} \left( \frac{qT^2}{\ell^{1/2}} + q^{3/4} T^{7/4} + T^{1+2\vartheta_0} (q^{7/8} \ell^{3/4} + q^{3/4} \ell) \right)
\end{split}
\end{displaymath}
by \eqref{M1} and \eqref{M2}. On the left hand side, we write
$$L^{(q)}(1/2 - i\tau, \bar{f}) = L(1/2 - i\tau, \bar{f})(1 - O(q^{-1/2}))$$
which introduces a total error of $O(q^{1/2+\varepsilon} T^{2+\varepsilon})$ by \eqref{largesieve}.
 
Now let $h$ be a $T$-admissible function as in Theorem \ref{thm3} with the additional property that it  has a   zero at $\pm i(1/2 \pm i\tau \pm n)$ for all sign combinations and all $n \in \Bbb{N}$, $1 \leq n \leq A+1$ and a quadruple zero at $\pm i(1/2 \pm i\tau)$. Then
$$h^{\ast}(t) :=  \frac{{\tt G}_{+ 1}(1/2 + i\tau, t) }{ {\tt G}_{+ 1}(1/2 -i\tau, t)}  h(t)$$
is a $T$-admissible function with a triple zero at $\pm i(1/2 \pm i\tau)$ so that the above reasoning is valid. Observing that 
$$ \frac{{\tt G}_{- 1}(1/2 -i\tau, t)}{{\tt G}_{- 1}(1/2 + i\tau, t)}\frac{{\tt G}_{+1 }(1/2 +i\tau, t)}{{\tt G}_{+1 }(1/2 - i\tau, t)}   
=   \frac{ 1 -   i \sinh(\pi \tau)/\cosh(\pi t)}{1 + i \sinh(\pi \tau)/\cosh(\pi t) } = Q_{\tau}(t) $$
completes the proof of Theorem \ref{thm3}. 

 \section{Amplification}
 
 We finally prove Theorem \ref{thm1} as a simple consequence of Theorem \ref{thm3}. Let $f_0\in \mathcal{B}(q, \chi)$ have spectral parameter $t$ and write $T = 1+|t|$. For any $f \in \mathcal{B}(q, \chi)$, we define 
$$A_f := \Bigl|\sum_{ \substack{p \leq L\\ p \nmid q} } \lambda_f(p) \bar{x}(p)\Bigr|^2 +  \Bigl|\sum_{\substack{p \leq L\\ p \nmid q}} \lambda_f(p^2) \bar{x}(p^2)\Bigr|^2, \quad x(n) = \text{sgn}(\lambda_{f_0}(n)) \in S^1 \cup \{0\}$$
for a parameter $L$ with $(qT)^{1/100} \leq L \leq (qT)^{1/2}$, say. 
Then  
\begin{equation}\label{lowerbound}
A_{f_0} = \Bigl(\sum_{\substack{ p\leq L\\ p \nmid q} } |\lambda_{f_0}(p) |\Bigr)^2 +  \Bigl(\sum_{\substack{ p\leq L\\ p \nmid q} } |\lambda_{f_0}(p^2)|  \Bigr)^2 \geq \frac{1}{2} \Bigl(\sum_{\substack{ p\leq L\\ p \nmid q} } |\lambda_{f_0}(p)| + |\lambda_{f_0}(p^2)| \Bigr)^2 \gg \frac{L^2}{\log L}
\end{equation}
by the prime number theorem and the Hecke relation $\lambda_{f_0}(p)^2 = \chi(p) + \lambda_{f_0}(p^2)$.  On the other hand,
\begin{displaymath}
\begin{split}
A_f =& \sum_{\substack{ p \leq L \\ p \nmid q}} (|x(p)|^2 + |x(p^2)|^2)+ \sum_{ \substack{p_1, p_2 \leq L \\ p_1p_2 \nmid q}} \left(\bar{x}(p_1)x(p_2)\bar{\chi}(p_2) + \delta_{p_1 = p_2} \bar{x}(p_1^2) x(p_2^2)\bar{\chi}(p_2)\right) \lambda_f(p_1p_2) \\
&+ \sum_{ \substack{p_1, p_2 \leq L \\ p_1p_2 \nmid q}}\bar{x}(p_1^2) x(p_2^2) \bar{\chi}(p_2^2)\lambda_f(p_1^2p_2^2).
\end{split}
\end{displaymath}
Employing the function  $h$ defined in \eqref{defh} (which we recall is positive for $t \in \Bbb{R} \cup [-i\vartheta_0, i \vartheta_0]$), we conclude from Theorem \ref{thm3} applied with $\ell \in \{1, p_1p_2, p_1^2p_2^2\}$ that
\begin{displaymath}
\begin{split}
 \sum_{ f \in \mathcal{B}(q, \chi) } & \begin{cases} 1, & \epsilon_f = 1\\ Q_{\tau}(t_f), & \epsilon_f = -1\end{cases} \Bigg\} A_f \frac{|L(1/2+  i\tau, f)|^4 } {L(1, \text{Ad}^2 f)}   h(t_f) \\
&\ll \left( qT^2L  + q^{1/2} T^2 L^2+ q^{3/4} T^{7/4} L^2+ T^{1+2\vartheta_0}q^{1/2}\left(q^{3/8} L^5 + q^{1/4} L^6 \right)\right)(TqL)^{\varepsilon}
\end{split}
\end{displaymath}
which implies in particular an upper bound for the real part and the imaginary part of the left hand side. 
Now we have for suitable constants $t_0 = t_0(\tau) > 0$, $t_1 = t_1(\tau) > 0$ that 
\begin{displaymath}
\begin{split}
 &   \Re Q_{\tau}(t) \geq 1/2, \quad t \in \Bbb{R}, |t| \geq t_0, \\
 &     \Re Q_{\tau}(t) \geq -1, \quad t \in \Bbb{R} \cup [-i\vartheta_0, i \vartheta_0],\\
 & -\text{sgn}(\tau) \cdot  \Im Q_{\tau}(t)  \geq t_1, \quad t \in [-t_0, t_0] \cup [-i\vartheta_0, i \vartheta_0].
 \end{split}
 \end{displaymath}
  Taking a suitable linear combination of real and imaginary part of the preceding inequality (namely real part minus $2 \text{ sgn}(\tau) t_1^{-1}$ times imaginary part), gives
  \begin{displaymath}
  \begin{split}
 & \sum_{ f \in \mathcal{B}(q, \chi) }   A_f \frac{|L(1/2+  i\tau, f)|^4 } {L(1, \text{Ad}^2 f)}   h(t_f) 
\\
&\ll \left( qT^2L  + q^{1/2} T^2 L^2+ q^{3/4} T^{7/4} L^2+ T^{1+2\vartheta_0}q^{1/2}\left(q^{3/8} L^5 + q^{1/4} L^6 \right)\right)(TqL)^{\varepsilon}. 
\end{split}
\end{displaymath}
Using  the  upper bound in \eqref{hl} and \eqref{lowerbound}, the left hand side is $\gg L^2 (qLT)^{-\varepsilon} |L(1/2+  i\tau, f_0)|^4 $ by positivity, and we arrive at
$$|L(1/2+  i\tau, f_0)|^4 \ll  \left( qT^2L^{-1} + q^{1/2} T^2+ q^{3/4} T^{7/4}  + T^{1+2\vartheta_0}q^{1/2}\left(q^{3/8} L^3 + q^{1/4} L^4 \right)\right)(TqL)^{\varepsilon}.$$
  We choose $L = q^{1/32} T^{(1-2\vartheta_0)/5}$ to conclude the proof of Theorem \ref{thm1}.

\end{document}